\documentclass[10pt]{thesis}

\usepackage{epsfig}
\usepackage{amsmath}
\usepackage{amssymb}
\usepackage{amsfonts}
\usepackage{amsthm}
\usepackage{latexsym}
\usepackage{mathrsfs}
\usepackage{amsxtra}

\usepackage[T1]{fontenc}
\usepackage{ae}

\input{xy}
\xyoption{all}

\newtheorem*{mainthmA}{Theorem A}
\newtheorem*{mainthmB}{Theorem B}

\newtheorem{theorem}{Theorem}[chapter]
\newtheorem{proposition}[theorem]{Proposition}
\newtheorem{lemma}[theorem]{Lemma}
\newtheorem{corollary}[theorem]{Corollary}
\newtheorem{conjecture}[theorem]{Conjecture}
\newtheorem{problem}[theorem]{Problem}

\theoremstyle{definition}
\newtheorem{definition}[theorem]{Definition}
\newtheorem{example}[theorem]{Example}

\theoremstyle{remark}
\newtheorem{remark}[theorem]{Remark}
\newtheorem{claim}{Claim}

\numberwithin{equation}{chapter}

\DeclareMathOperator{\Diff}{Diff}

\DeclareMathOperator{\M}{\mathfrak{M}}

\DeclareMathOperator{\Per}{Per}

\DeclareMathOperator{\SL}{SL}
\DeclareMathOperator{\GL}{GL}

\DeclareMathOperator{\Dis}{\mathcal{D}^\prime}
\DeclareMathOperator{\rank}{rank}

\DeclareMathOperator{\diver}{div}
\DeclareMathOperator{\cl}{cl}
\DeclareMathOperator{\sign}{sign}

\newcommand{\Lie}{\mathcal{L}}
\newcommand{\abs}[1]{\lvert {#1}\rvert}
\newcommand{\bb}[1]{\mathbb{#1}}
\newcommand{\scr}[1]{\mathscr{#1}}
\newcommand{\X}{\mathfrak{X}}

\def\ie{i.e.\ }

\def\etc{\textit{etc.}}

\title{Toward the classification of cohomology-free vector fields}
\author{Alejandro Kocsard}

\submitted{Rio de Janeiro, 23 de maio de 2007}

\abstract{Given a smooth vector field $X$ on a closed orientable
  $d$-manifold $M$, many questions about the dynamics of its induced
  flow can be studied analyzing the following \emph{cohomological
    equation}:
  \begin{displaymath}
    \Lie_X u = \xi,
  \end{displaymath}
  where $\xi$ is a given real function on $M$, $u\colon M\to\bb{R}$ is
  the solution that we look for (in a certain regularity class) and
  $\Lie_X$ is the Lie derivative in the $X$ direction.

  In 1984, Anatole
  Katok~\cite{rigid-problem-list,katok-robinson,katok} proposed to
  characterize those vector fields which are cohomologicaly trivial.
  More precisely, he conjectured that if $X$ is so that for all smooth
  function $\xi\colon M\to\bb{R}$, there exist a constant
  $c=c(\xi)\in\bb{R}$ and $u\in C^\infty(M,\bb{R})$ verifying
  \begin{displaymath}
    \Lie_X u = \xi - c,
  \end{displaymath}
  then $X$ should be smoothly conjugated to a Diophantine (constant)
  vector field on $\bb{T}^d$. In particular, $M$ should be
  diffeomorphic to $\bb{T}^d$.

  The main goal of this work is to prove the validity of Katok
  Conjecture for $3$-manifolds.}

\acknowledgements{Aos meus ``viejos'', pelo carinho e ap\'oio
  incondicional.

  \`A Marina, ``pelo simples coragem de me querer''.

  Ao IMPA como institui\c{c}\~{a}o, e em particular ao Welington de
  Melo, pela excelente orienta\c{c}\~{a}o; ao Enrique Pujals, pela sua
  amizade, paci\^{e}ncia e aten\c{c}\~{a}o que tem me dedicado
  durantes esses anos; ao Jacob Palis, pela sua confian\c{c}a e
  ap\'oio; ao Marcelo Viana, pela sua constante preocupa\c{c}\~{a}o
  com a nossa evolu\c{c}\~{a}o.  Tamb\'{e}m desejo agradecer aos meus
  colegas de turma, Meysam, Jimmy, Martin, Paulo, e muito
  espacialmente ao Andr\'es Koropecki, pela sua amizade e por tudo o
  que temos compartilhado ao longo desses quatro anos no Rio.

  Aos membros da banca, Viviane Baladi, Carlos Gustavo Moreira (Gugu),
  e muito especialmente a Mar\'{\i}a Alejadra Rodrigu\'{e}z Hertz
  (Jana), pelas muitas conversas esclarecedoras sobre a Conjectura de
  Katok.

  Tamb\'{e}m gostaria de agradecer ao Giovanni Forni pelas suas muitas
  sugest\~{o}es e disponibilidade para discutir sobre os assuntos
  apresentados neste trabalho.

  Aos meus amigos ``brasileiros'', Rudy, Gabriel, Javier, Hern\'{a}n,
  Juan, Vero, Ro, Jose, Luis, Tati, Alexandre, Ximena, Mercedes,
  Katia, Ramiro, Diana, Andr\'{e}, Ernesto, Fernanda, Isabella e Gil.

  Aos meus amigos da matem\'{a}tica de Rosario, Santiago, Guillermina,
  Lisandro (e Gabi), Mariela (e Juanca) Cristian, Pini e \'{E}rica.

  Ao CNPq e \`{a} FAPERJ, pelo ap\'{o}io econ\'{o}mico recebido.

  Aos meus professores Pedro Marangunic, Nidia Jeifez, Rafael Verdes e
  Hugo \\ Aimar por terem tido um papel central na minha
  forma\c{c}\~{a}o durante a Licenciatura en Matem\'{a}tica na UNR.

  Finalmente, gostaria de agradecer \`{a} Olimpiada Matem\'{a}tica
  Argentina (que injustamente nunca aperece dentro da se\c{c}\~{a}o
  ``Forma\c{c}\~{a}o Acad\'{e}mica'' do meu CV), especialmente a
  \'Erica Hinrichsen, Pablo Lotito, Patricia Fauring e Flora
  Guti\'errez, por ter sido a minha primeira escola de
  matem\'{a}tica.}

\begin{document}

\chapter{Introduction}
\label{chap:intro}

The main goal in Differentiable Dynamics consists in understanding the
global behavior of ``most'' of the orbits of systems, where the
\emph{phase space} is represented by a compact differential manifold
$M$, and the evolution by a diffeomorphism $f\in\Diff^r(M)$ (discrete
time case), or by a $C^r$ flow $\Phi\colon M\times\bb{R}\to M$
(continuous time case).

Looking for the unification of the notation, we can assume that the
dynamics of our system is given by a $C^r$ Lie group action on $M$.
More precisely, if $(\bb{G},+)$ denotes any analytic abelian Lie
group, we shall suppose that $\bb{G}$ represents the time and the
evolution of the system is given by a $C^r$ $\bb{G}$-action
$\Gamma\colon M\times \bb{G}\to M$.

\section{Cocycles and Coboundaries}
\label{sec:cocycles-cobound}

When we analyze different questions about the dynamics of $\Gamma$,
there is a family of objects that appears repeatedly and in a very
natural way (see Section~\ref{sec:cohomo-eq} for some examples). These
are the \emph{real cocycles over} $\Gamma$:

\begin{definition}
  \label{def:cocycles}
  Given a $C^r$ $\bb{G}$-action $\Gamma\colon M\times\bb{G}\to M$, a
  \emph{real cocycle over} $\Gamma$ (or simply a \emph{cocycle}) is a
  $C^k$ map (usually $k\leq r$) $\Xi\colon M\times\bb{G}\to\bb{R}$
  such that
  \begin{equation}
    \label{eq:def-cocycles}
    \Xi(p,g_0+g_1)=\Xi(\Gamma(p,g_0),g_1)+\Xi(p,g_0),\quad\forall p\in
    M,\ \forall g_0,g_1\in\bb{G}. 
  \end{equation}
\end{definition}

Within this framework it is natural to consider the following
equivalence relation between cocycles:

\begin{definition}
  \label{def:cocycles-cohomolog}
  We shall say that two $C^k$ cocycles $\Xi,\Theta\colon
  M\times\bb{G}\to M$ over $\Gamma$ are $C^s$-\emph{cohomologous}
  (with $s\leq k\leq r$) if there exists a $C^s$ map $\alpha\colon
  M\to\bb{R}$ verifying
  \begin{displaymath}
    \Xi(p,g)=\alpha(\Gamma(p,g))+\Theta(p,g)-\alpha(p),\quad\forall
    p\in M,\ \forall g\in\bb{G}.
  \end{displaymath}
\end{definition}

On the other hand, notice that given any real $C^k$ function
$\xi\colon M\to\bb{R}$, we can easily construct a cocycle $\Xi$ over
$\Gamma$ defining
\begin{equation}
  \label{eq:def-coboundaries}
  \Xi(p,g)\doteq \xi(\Gamma(p,g)) - \xi(p),\quad\forall p\in M,\
  \forall g\in\bb{G}.
\end{equation}

Cocycles constructed as above are very important and deserve a special
name: they are called \emph{coboundaries}.  In other words, we may say
that a cocycle is a coboundary if and only if it is cohomologous to
the null cocycle.

These names come from (abstract) Group Cohomology Theory. In fact, if
we suppose that $\Gamma$ is $C^\infty$, then it induces in a natural
way a $\bb{G}$-action on $C^\infty(M,\bb{R})$, turning
$C^\infty(M,\bb{R})$ into a $\bb{G}$-module. So, in a purely algebraic
way we can define the cohomology complex $H^*(\bb{G},\Gamma)$ (see
\cite{atiyah} for example). In this way, $H^1(\bb{G},\Gamma)$ happens
to be canonically isomorphic to the quotient vector space of all
smooth cocycles over $\Gamma$ by the subspace of all smooth
coboundaries. However, since in the future we shall not make any other
reference to higher cohomology groups, the reader can simply consider
this algebraic construction as a justification for the chosen names.

\section{Cohomological Equations}
\label{sec:cohomo-eq}

Since we are mainly interested in the ``classical group actions'',
from now on we shall assume that $\Gamma$ is a differentiable
$\bb{G}$-action, being $\bb{G}=\bb{Z}$ or $\bb{R}$.

As it was already mentioned, cocycles appear naturally in different
contexts when we want to study some dynamical properties of $\Gamma$.
Among the problems in Differentiable Dynamics that can be reduced to
cohomological considerations we can mention:
\begin{enumerate}
\item \emph{Existence of invariant volume forms} (see Section 5.1 in
  the book of Katok and Hasselblat \cite{katok-hassel}).
\item \emph{Stability of hyperbolic torus automorphisms} (see Section
  2.6 in \cite{katok-hassel}).
\item \emph{Liv\v sic Theory} (see Section 19.2 in
  \cite{katok-hassel}, Section 3.4 in the survey of Katok and Robinson
  \cite{katok-robinson}, or the original work Liv\v sic \cite{liv}).
\item \emph{KAM Theory} (see the survey of R. de la Llave
  \cite{delallave}).
\item \emph{Constructions of minimal conservative but non uniquely
    ergodic diffeomorphisms} (see the classical work of H. Furstenberg
  \cite{furs}).
\end{enumerate}

In all the cases listed above the main problem consists in proving
that a given cocycle is or is not a coboundary, or more generally,
that it is or it is not $C^s$-cohomologous to another given cocycle.

This is the reason why it is so important to analyze the existence of
solutions $u\colon M\to\bb{R}$ (in a certain regularity class) for the
following difference equation:
\begin{equation}
  \label{eq:gen-cohomo-eq}
  u(\Gamma(p,g))-u(p)=\Xi(p,g),\quad\forall p\in M,\ \forall
  g\in\bb{G},
\end{equation}
where $\Xi$ is a given real cocycle over the $\bb{G}$-action $\Gamma$.
These equations deserve a special name:
\begin{definition}
  A difference equation like (\ref{eq:gen-cohomo-eq}) will be called a
  \emph{cohomological equation}.
\end{definition}

In the particular case that $\bb{G}=\bb{Z}$, the cocycle $\Xi$ is
``generated'' by the function $\xi(p)\doteq\Xi(p,1)$. Indeed, it holds
\begin{displaymath}
  \Xi(p,n)=
  \begin{cases}
    0,\text{ if } n=0, \\
    \sum_{i=0}^{n-1}\xi(f^i(p)), \text{ if } n>0, \\
    -\sum_{i=n}^{-1}\xi(f^i(p)), \text{ if } n<0,
  \end{cases}
\end{displaymath}
where $f\doteq\Gamma(\cdot,1)\in\Diff^r(M)$. And so, in this case the
cohomological equation~(\ref{eq:gen-cohomo-eq}) can be written as
\begin{equation}
  \label{eq:Z-cohomo-eq}
  u\circ f - u=\xi.
\end{equation}

On the other hand, when $\bb{G}=\bb{R}$ the cocycle $\Xi$ has an
``infinitesimal generator'' defined by
\begin{displaymath}
  \xi(p)\doteq \partial_t \Xi(p,t)\Big|_{t=0},
  \quad\forall p\in M.
\end{displaymath}
In this case, $\Xi(p,t)=\int_0^t\xi(\Gamma(p,s))\:\mathrm{d}s$, and
hence, derivating equation~(\ref{eq:gen-cohomo-eq}) with respect to
the time variable, we get the following differential equation:
\begin{equation}
  \label{eq:R-cohomo-eq}
  \Lie_X u = \xi,
\end{equation}
where $X\in\X(M)$ is the vector field generating $\Gamma$, \ie
$X(p)\doteq\partial_t\Gamma(p,t)\big|_{t=0}$, and $\Lie_X$ denotes the
Lie derivative along $X$.

\section{Obstructions}
\label{sec:obstruct}

In general it is not an easy task to determine if a particular
cohomological equation admits some solution in a particular regularity
class. So, it appears as an important problem to characterize the
\emph{``set of obstructions''} for the existence of ($L^p$, $C^r$,
\etc ) solutions for equations like (\ref{eq:gen-cohomo-eq}).

For example, if $\Gamma$ is a given $\bb{R}$-action, then the very
first obstructions that we can find for the existence of continuous
solutions for equation~(\ref{eq:R-cohomo-eq}) are the elements of
$\M(\Gamma)$, the set of Borel finite measures on $M$ which are left
invariant by the flow $\Gamma$.  More precisely, if $u$ is a
continuous solution of equation~(\ref{eq:R-cohomo-eq}), since
\begin{equation*}
  \frac{1}{T}\int_0^T\xi(\Gamma(p,t))\:\mathrm{d}t=
  \frac{1}{T} \big(u(\Gamma(p,t))-u(p)\big)\xrightarrow{T\to\infty}0,
\end{equation*}
as a straight-forward consequence of Birkhoff ergodic theorem we have
that
\begin{equation*}
  \int_M\xi\:\mathrm{d}\mu=0,\quad\text{for every }\mu\in\M(\Gamma).
\end{equation*}

In Section~\ref{sec:stri-erg} we shall see that the set of
$\Gamma$-invariant distributions, in the sense of Schwartz, is the
most natural space for looking for obstructions for the existence of
smooth solutions for equation~(\ref{eq:R-cohomo-eq}) (or
(\ref{eq:Z-cohomo-eq})).

Two very classical results which completely characterize this
\emph{set of obstructions} in two particular, and in some sense,
extremal opposite situations, are due to Gottschalk and
Hedlund~\cite{got-hed} and to Liv\v sic~\cite{liv}.

In the first one, if $\Gamma$ is a continuous minimal $\bb{Z}$-action
(\ie every point in $M$ has a dense $\Gamma$-orbit) generated by a
homeomorphism $f$ on $M$ and if $\xi\in C^0(M,\bb{R})$, Gottschalk and
Hedlund proved that equation~(\ref{eq:Z-cohomo-eq}) admits a
continuous solution $u$ if and only if the family of functions

\begin{displaymath}
  \bigg\{\sum_{i=0}^{n-1}\xi\circ f^i\bigg\}_{n\geq 1}
\end{displaymath}
is uniformly bounded in $C^0(M,\bb{R})$.

In the second one, Liv\v sic studied the case where $\Gamma$ is a
$C^2$ hyperbolic $\bb{R}$-action (\ie an Anosov flow) induced by a
vector field $X\in\X^3(M)$. Assuming that $\xi$ is H\"older
continuous, he proved that the only obstruction for the existence of a
H\"older continuous solution $u$ for equation~(\ref{eq:R-cohomo-eq})
is given by the set of probabilities concentrated on the periodic
orbits, \ie there is a H\"older continuous solution $u$ as long as

\begin{displaymath}
  \int_0^{\tau(z)}\psi(\Gamma(z,s))\:\mathrm{d}s=0,\quad\forall
  z\in\Per(\Gamma), 
\end{displaymath}
and where $\tau(z)\doteq\inf\{t>0 : \Gamma(z,t)=z\}$.

It is interesting to remark that both results \cite{got-hed} and
\cite{liv} hold for $\bb{R}$-actions as well as for $\bb{Z}$-actions.

There are more recent results that completely characterize the sets of
obstructions in some other cases. For example, cohomological equations
associated to area-preserving flows on higher genus surfaces have been
studied by Giovanni Forni~\cite{forni1,forni2} (the torus case is
rather classical); and associated to their ``very close relatives'',
the interval exchanged maps, by Stefano Marmi, Pierre Moussa and
Jean-Christophe Yoccoz~\cite{mmy1,mmy2}. Other very important flows
that nowadays are very well understood from the cohomological point of
view are homogeneous ones on nilmanifolds: this study is due to Livio
Flaminio and Giovanni Forni. They started studying some particular
cases (horocycle flows, nilflows on Heisenberg manifolds) in
\cite{fla-for1,fla-for2}, and the general case was settled in
\cite{fla-for}.

\section{Cohomology-free Dynamical Systems}
\label{sec:cohomo-free-vect}

As it was already mentioned in Section~\ref{sec:obstruct}, in general
it is very difficult to characterize the set of obstructions for the
existence of solutions for a cohomological equation like
(\ref{eq:R-cohomo-eq}). With the aim of understanding the nature
(topological, analytical, \etc) of this set of obstructions, Anatole
Katok, in the early `80, proposed the following

\begin{definition}
  \label{def:cohomo-free-def}
  Given a closed manifold $M$, we say that a smooth $\bb{G}$-action
  $\Gamma\colon M\times\bb{G}\to\bb{M}$ is \emph{cohomology-free} if
  any smooth real cocycle over $\Gamma$ is $C^\infty$-cohomologous to
  a constant one.
\end{definition}

Notice that two different constant cocycles are never smoothly
cohomologous, and so, the first cohomology group of any smooth action
always contains a subgroup isomorphic to $\bb{R}$. Therefore, we can
say that a smooth action is cohomology-free if and only if its first
cohomology group is as small as possible.

For the sake of clarity of the exposition, from now on we shall mainly
concentrate on smooth $\bb{R}$-actions, \ie flows induced by
$C^\infty$ vector fields. In this particular case,
Definition~\ref{def:cohomo-free-def} can be restated in the following
way:

We say that $X\in\X(M)$ is \emph{cohomology-free} if given any $\xi\in
C^\infty(M,\bb{R})$, there exist a constant $c(\xi)\in\bb{R}$ and
$u\in C^\infty(M,\bb{R})$ verifying
\begin{equation}
  \label{eq:cohomology-free-def}
  \Lie_X u = \xi - c(\xi).
\end{equation}

\begin{remark}
  \label{rem:C-conj}
  It is clear that the set of cohomology-free vector fields is closed
  under $C^\infty$-conjugacy.
\end{remark}

To introduce the prototypical example of cohomology-free vector
fields, first we need to state the following

\begin{definition}
  \label{def:diophantine}
  We say that $\alpha=(\alpha_1,\ldots,\alpha_d)\in\bb{R}^d$ is a
  \emph{Diophantine vector} if there exist real constants $C,\tau>0$
  satisfying
  \begin{equation}
    \label{eq:diophantine-def}
    \bigg|\sum_{i=1}^d \alpha_i p_i\bigg| > C \bigg(\max_{1\leq i\leq d}
    \abs{p_i}\bigg)^{-\tau}, 
  \end{equation}
  for every $p=(p_1,\ldots,p_d)\in\bb{Z}^d\setminus\{0\}$.

  A vector field $X_\alpha$ on the $d$-dimensional torus $\bb{T}^d$
  verifying $X_\alpha\equiv\alpha$ will be called a \emph{Diophantine
    vector field}.
\end{definition}

\begin{example}
  \label{ex:diop-vect}
  Diophantine vector fields on tori are cohomology-free.

  In fact, let $\alpha\in\bb{R}^d$ be a Diophantine vector. The Haar
  measure on $\bb{T}^d$ is the only $X_\alpha$-invariant probability
  measure and if $\xi\colon\bb{T}^d\to\bb{R}$ is an arbitrary
  $C^\infty$ function, considering its Fourier expansion
  \begin{displaymath}
    \xi(\theta)=\sum_{k\in\bb{Z}^d}\hat{\xi}_k e^{2\pi i k\cdot\theta},
  \end{displaymath}
  we can define $u$, at first just formally, writing
  \begin{displaymath}
    u(\theta)=\sum_{k\in\bb{Z}^n\setminus\{0\}}
    \frac{\hat{\xi}_k}{k\cdot\alpha}e^{2\pi ik\cdot\theta}.
  \end{displaymath}
  Taking into account estimate~(\ref{eq:diophantine-def}), we easily
  see that $u\in C^\infty(\bb{T}^n,\bb{R})$ and, by construction, it
  holds
  \begin{displaymath}
    \Lie_{X_\alpha} u = \xi - \hat{\xi}_0.
  \end{displaymath}
\end{example}

As we will see in Theorem~\ref{thm:cf-vect-tori}, these are the only
cohomology-free vector fields on tori, of course, modulo
$C^\infty$-conjugacy.

Considering this example and the previous work of Stephen Greenfield
and Nolan Wallach \cite{hypo-vect} on globally hypoelliptic vector
fields (see Section~\ref{sec:glob-hypo-vect} for precise definitions),
Anatole Katok proposed in \cite{rigid-problem-list} the following
conjecture characterizing the cohomology-free vector fields:

\begin{conjecture}[Katok Conjecture \cite{rigid-problem-list}]
  \label{conj:katok}
  If $M$ is a compact, connected, orientable $d$-manifold, and $X$ is
  a cohomology-free vector field on $M$, then $M$ is diffeomorphic to
  the torus $\bb{T}^d$, and therefore, $X$ is $C^\infty$ conjugated to
  a Diophantine constant vector field on $\bb{T}^d$.
\end{conjecture}

Some results supporting Katok Conjecture have recently appeared.
First, Federico and Jana Rodriguez-Hertz \cite{fede-jana} have proved
that a manifold supporting a cohomology-free vector field must fiber
over the torus of dimension equal to the first Betti number of the
manifold (see Theorem~\ref{thm:fede-jana} for the precise statement);
and secondly, Livio Flaminio and Giovanni Forni \cite{fla-for} have
proved that tori are the only nilmanifolds supporting cohomology-free
homogeneous flows.

\section{Main Results and Outline of this Work }
\label{sec:out-main-results}

The main goal of this work is to present a complete proof of Katok
Conjecture in dimension $3$. For this, the rest of the work will be
organized as follows:

In Chapter~\ref{chap:gen} we shall present general properties about
cohomology-free vector fields, some of which are very classical, like
strict ergodicity, and the rather new result due to Federico and Jana
Rodr\'\i guez-Hertz, Theorem~\ref{thm:fede-jana}.

In Chapter~\ref{chap:b_1>0} we shall expose the proof of our first
result toward the classification of cohomology-free vector fields on
$3$-manifolds:

\begin{mainthmA}
  Let $M$ be a closed and orientable $3$-manifold verifying
  \begin{equation*}
    \beta_1(M)=\dim H_1(M,\bb{Q})\geq 1,
  \end{equation*}
  and suppose that there exists a smooth cohomology-free vector field
  $X\in\X(M)$.  Then $M$ is diffeomorphic to $\bb{T}^3$ and $X$ is
  $C^\infty$-conjugated to a Diophantine constant vector field.
\end{mainthmA}

While this work was in progress, Giovanni Forni~\cite{forni-weins}
communicated to us that he had proved the following
result\footnote{Forni independently got a proof of Theorem A, too.}:

\begin{mainthmB}
  \label{thm:forni-0-betti}
  If $M$ is a closed orientable $3$-manifold with $H_1(M,\bb{Q)}=0$
  and $X\in\X(M)$ is a cohomology-free vector field, then there exists
  a $1$-form $\alpha$ on $M$ verifying
  \begin{displaymath}
    \alpha\wedge d\alpha\neq 0,\quad i_X\alpha\equiv 1,\quad i_X
    d\alpha\equiv 0.
  \end{displaymath}
  In other words, $\alpha$ is a contact form and $X$ is its induced
  Reeb vector field.
\end{mainthmB}

On the other hand, Clifford Taubes~\cite{taubes} has recently proved
Weinstein Conjecture which asserts that every Reeb vector field on a
$3$-manifold must exhibit a periodic orbit. This clearly contradicts
the minimality (see Corollary~\ref{cor:minimality}) of the flow
induced by a cohomology-free vector field. Therefore, Theorem~B lets
us affirm that there is no cohomology-free vector field on
$3$-manifolds with vanishing first Betti number and thus we get

\begin{corollary}[Katok Conjecture in dimension $3$]
  If $M$ is closed and orientable $3$-manifold and $X\in\X(M)$ is a
  smooth cohomology-free vector field on $M$, then $M$ is
  diffeomorphic to $\bb{T}^3$ and $X$ is $C^\infty$-conjugated to a
  constant Diophantine vector field.
\end{corollary}

In Chapter~\ref{chap:b_1=0} we will sketch very briefly Forni's proof
of Theorem~B (that he kindly communicated to us) and we will present
another proof, using completely different techniques. We hope this can
help to get a better comprehension of the whole problem.

For the sake of completeness, in that chapter we will also recall some
fundamental facts about Contact Geometry and we shall precisely state
Weinstein Conjecture.

Finally, in Chapter~\ref{chap:remarks-problems} we propose some open
problems and consider some final remarks about the results presented
in this work.

\section{Notation and Conventions}
\label{sec:not-conv}

For simplicity, we will mainly work in the $C^\infty$ category and the
word \emph{smooth} will be used as a synonymous of $C^\infty$.

We shall say that a manifold is \emph{closed} if it is compact,
connected and its boundary is empty.

Along this work, $M$ will denote a smooth closed orientable
$d$-dimensional manifold, and most of the time $d=3$.

The linear space of all $C^r$ vector fields on $M$ will be denoted by
$\X^r(M)$, and to simplify the notation, we shall just write $\X(M)$
for the space of smooth vector fields.

Analogously, $\Diff^r(M)$ will stand for the set of $C^r$
diffeomorphism of $M$ and we will simply write $\Diff(M)$ for the set
of smooth diffeomorphisms.

The expression $\Lambda^k(M)$ will be used for the space of smooth
$k$-forms on $M$, and given any $X\in\X(M)$,
$i_X\colon\Lambda^k(M)\to\Lambda^{k-1}(M)$ shall denote the
contraction by $X$ (also called interior product).

As usual, we shall identify $\Lambda^0(M)$ with $C^\infty(M,\bb{R})$.

Given any $X\in\X(M)$, $\{\Phi_X^t\}_{t\in\bb{R}}$ will denote the
flow induced by $X$.

If $T$ denotes any smooth tensor field on $M$, the Lie derivative of
$T$ along $X$ will be denoted by $\Lie_X T$ and defined by
\begin{displaymath}
  \Lie_X T(x)\doteq \lim_{t\to 0}\frac{(\Phi_X^t)^*T(x)-T(x)}{t},
  \quad\forall x\in M.
\end{displaymath}

The set of all finite signed Borel measures on $M$ (\ie real
continuous linear functionals on $C^0(M,\bb{R})$) shall be denoted by
$\M(M)$, and we will write $\Dis(M)$ for the space of all real
continuous linear functionals on $C^\infty(M,\bb{R})$.

Given any smooth fibration $p\colon N \to M $, the fiber over any
$x\in M$ shall be denoted by $N_x$, and we will write $\Gamma(N)$ for
the space of smooth sections (\ie maps $s\colon M\to N$ verifying
$p\circ s=id_M$). The only exception for this notational convention is
the tangent bundle over $M$: in this case $\pi\colon TM\to M$ will
denote the canonical projection and we will write $T_xM$ for
$\pi^{-1}(x)$ and $\X(M)$ for $\Gamma(TM)$.

Similarly, given any foliation $\scr{F}$ on $M$, $\scr{F}_x$ or
$\scr{F}(x)$ will stand for the leaf of $\scr{F}$ through $x\in M$.

It is very important to remark that, in order to avoid confusions
along this work we shall use the term \emph{distribution} in the
``sense of Schwartz,'' \ie for us a \emph{distribution} will be any
element of $\Dis(M)$. This word has a completely different meaning in
Differential Geometry. Indeed, we shall use the expression
\emph{k-plane field}, or \emph{line field} when $k=1$, to denote the
objects that are commonly named distributions in Differential
Geometry.

There are some relationships between the linear spaces
$C^\infty(M,\bb{R})$, $\Lambda^d(M)$, $\M(M)$ and $\Dis(M)$.  First,
since the elements of $\M(M)$ can be considered as linear continuous
functionals on $C^0(M,\bb{R})$, it can be canonically embedded in
$\Dis(M)$. On the other hand, it is very easy to see that each element
of $\Lambda^d(M)$ (where $d=\dim M$) naturally induces a signed
measure, \ie we can assume that $\Lambda^d(M)\subset\M(M)$. And
finally, since $M$ is supposed to be orientable, we can choose a
volume form on $M$ and use it to get a bijection between
$C^\infty(M,\bb{R})$ and $\Lambda^d(M)$.  However, it is important to
remark that in this case this identification is not canonical at all.

Since we have defined the Lie derivative $\Lie_X$ on
$C^\infty(M,\bb{R})$, we can easily extend it by duality to $\Dis(M)$.
In fact, we can define $\Lie_X\colon\Dis(M)\to\Dis(M)$ writing
\begin{displaymath}
  \langle \Lie_X T,\psi \rangle \doteq - \langle
  T,\Lie_X\psi\rangle, \quad\forall T\in\Dis(M),\ \forall\psi\in
  C^\infty(M,\bb{R}). 
\end{displaymath}

In this way, it is reasonable to define the set of $X$-invariant
distributions and measures by

\begin{align*}
  \Dis(X) &\doteq \{T\in\Dis(M):\Lie_X T=0\}, \\
  \M(X) & \doteq \{\mu\in\M(M):(\Phi_X^t)_*\mu=\mu,\ \forall
  t\in\bb{R}\} \\
  &=\{\mu\in\M(M): \mu\in\Dis(X)\}
\end{align*}

Finally, the $d$-dimensional torus will be denoted by $\bb{T}^d$ and
the quotient Lie group $\bb{R}^d/\bb{Z}^d$ will be our favorite model
for it. $\mathrm{pr}_{\bb{Z}^d}\colon\bb{R}^d\to\bb{T}^d $ will denote
the canonical quotient projection. The Haar probability measure on
$\bb{T}^d$, also called the Lebesgue measure, will be denoted by
$\mathrm{Leb}^d$.

In general, an arbitrary point of $\bb{T}^d$ shall be denoted by
$\theta=(\theta^0,\theta^1,\ldots,\theta^{d-1})$.

It is a very well-known fact that there exists a canonical group
isomorphism between the group of automorphisms of $\bb{T}^d$ and
$\GL(d,\bb{Z})$. Taking this into account, if $A\in\Diff(\bb{T}^d)$ is
any Lie group automorphism of $\bb{T}^d$, the corresponding element of
$\GL(d,\bb{Z})$ will be denoted by $\hat{A}$.  Notice that $A$ and
$\hat{A}$ are related by $A\circ\mathrm{pr}_{\bb{Z}^d}=
\mathrm{pr}_{\bb{Z}^d}\circ\hat{A}$.


\chapter{General Properties of Cohomology-free Vector Fields}
\label{chap:gen}

\section{Strict Ergodicity}
\label{sec:stri-erg}

This section is devoted to proving that every cohomology-free vector
field is strictly ergodic, \ie it is uniquely ergodic and every orbit
of its induced flow is dense on the whole manifold.  These results are
very classical and, as the reader will see, the proofs are rather
simple. Nevertheless, we decided to include them here for the sake of
completeness and with the purpose of making easier the reading of this
work.

As it was already mentioned in Section~\ref{sec:not-conv}, we shall
assume that $M$ is a closed orientable $d$-manifold.

\begin{proposition}
  \label{pro:uniq-ergodic}
  If $X\in\X(M)$ is a cohomology-free vector field, then its induced
  flow $\{\Phi_X^t\}$ is uniquely ergodic.
\end{proposition}

\begin{proof}
  Let $\psi\colon M\to\bb{R}$ be any smooth function and let
  $c(\psi)\in\bb{R}$ and $u\in C^\infty(M,\bb{R})$ be as in
  equation~(\ref{eq:cohomology-free-def}). Then we have
  \begin{equation}
    \label{eq:bir-sum}
    \frac{1}{T}\bigg( \int_0^T\psi(\Phi_X^s(p))\:\mathrm{d}s \bigg)
    = \frac{1}{T}\left(u(\Phi_X^T(p))-u(p)\right) + c(\psi), 
  \end{equation}
  for every $p\in M$ and every $T>0$.

  Then, if $\mu$ is an arbitrary $X$-invariant ergodic probability
  measure, by the Birkhoff ergodic theorem we know that the left side
  of equation~(\ref{eq:bir-sum}) must converge to
  $\int\psi\:\mathrm{d}\mu$, for $\mu$-almost every $p\in M$, when
  $T\rightarrow\infty$. On the other hand, since $u$ is bounded, the
  right side of~(\ref{eq:bir-sum}) converges to $c$. Therefore,
  $\int\psi\:\mathrm{d}\mu=c(\psi)$, for every $\mu\in\M(X)$, and
  since $C^\infty(M,\bb{R})$ is dense in $C^0(M,\bb{R})$, we conclude
  that $\M(X)$ contains only one element.
\end{proof}

In fact, we can prove a stronger result:

\begin{proposition}
  \label{pro:one-dim-invar-dist}
  If $X\in\X(M)$ is a cohomology-free vector field, then
  \begin{displaymath}
    \dim \Dis(X) = 1.
  \end{displaymath}
\end{proposition}

\begin{proof}
  Given an arbitrary $\psi\in C^\infty(M,\bb{R})$, let $u\in
  C^\infty(M,\bb{R})$ and $c(\psi)\in\bb{R}$ be such that
  \begin{displaymath}
    \Lie_X u = \psi - c(\psi).
  \end{displaymath}

  Then, for any $T\in\Dis(X)$ we have
  \begin{displaymath}
    \langle T,\psi \rangle = \langle T,\Lie_X u + c(\psi)\rangle =
    -\langle \Lie_X T, u\rangle + \langle T,c(\psi)\rangle = \langle
    T,c(\psi)\rangle.
  \end{displaymath}
  From this we can easily conclude that $\dim\Dis(X)=1$.
\end{proof}

We can also get the following regularity result for the elements of
$\Dis(X)$:

\begin{proposition}
  \label{pro:invar-volume-form}
  Let $X\in\X(M)$ be a cohomology-free vector field. Then there exists
  a smooth volume form $\Omega\in\Lambda^d(M)$ such that
  $\Lie_X\Omega\equiv 0$.
\end{proposition}

\begin{proof}
  Since we are assuming that $M$ is orientable, let
  $\tilde{\Omega}\in\Lambda^d(M)$ be an arbitrary smooth volume form.
  Let us define $\diver_{\tilde{\Omega}}X\in C^\infty(M,\bb{R})$ as
  the only smooth function verifying
  \begin{displaymath}
    \Lie_X\tilde{\Omega}=(\diver_{\tilde{\Omega}}X)\tilde{\Omega}.
  \end{displaymath}
  
  Hence there exist a smooth function $u\colon M\to\bb{R}$ and a real
  constant $c=c(\diver_{\tilde{\Omega}}X)$ satisfying
  \begin{displaymath}
    \Lie_X u = -(\diver_{\tilde{\Omega}}X) + c.
  \end{displaymath}

  Therefore, if we define $\Omega\doteq\exp(u)\tilde{\Omega}$, we
  obtain
  \begin{equation*}
    \begin{split}
      \Lie_X\Omega & =e^u\Lie_X\tilde{\Omega} +
      (e^u\Lie_Xu)\tilde{\Omega} \\
      &=e^u(\diver_{\tilde{\Omega}}X)\tilde{\Omega} +
      e^u(-(\diver_{\tilde{\Omega}}X) + c)\tilde{\Omega} \\
      & = ce^u\tilde{\Omega}=c\Omega.
    \end{split}
  \end{equation*}

  Finally, this clearly implies that
  $(\Phi_X^t)^*\Omega=(1+tc)\Omega$, and since the total
  $\Omega$-volume of $M$ is invariant, we have $c=0$.
\end{proof}

As a direct consequence of Propositions~\ref{pro:uniq-ergodic}
and~\ref{pro:invar-volume-form} we get the following

\begin{corollary}
  \label{cor:minimality}
  If $X\in\X(M)$ is a cohomology-free vector field, then the induced
  flow $\{\Phi_X^t\}_{t\in\bb{R}}$ is minimal, \ie it holds
  \begin{displaymath}
    \cl\big\{\Phi_X^t(p) : t\in\bb{R}\big\}=M,\quad\forall p\in M.
  \end{displaymath}
\end{corollary}

\section{Cohomology-free Vector Fields on Tori}
\label{sec:cf-vect-tori}

The aim of this section consists in proving that the Diophantine
vector fields are the only cohomology-free ones on tori, modulo
$C^\infty$-conjugacy.  More precisely, we shall prove the following

\begin{theorem}
  \label{thm:cf-vect-tori}
  If $X\in\X(\bb{T}^d)$ is a cohomology-free vector field on
  $\bb{T}^d$, then there exist a Diophantine vector $\alpha\in\bb{R}$
  (see Definition~\ref{def:diophantine}) and $f\in\Diff(\bb{T}^d)$
  homotopic to the identity such that
  \begin{displaymath}
    Df(X(\theta))\equiv\alpha.
  \end{displaymath}
\end{theorem}

This result is essentially due to Richard Luz and Nathan dos Santos.
In fact, in \cite{luz-santos} they proved that the only
cohomology-free diffeomorphisms on $\bb{T}^d$ homotopic to the
identity are those $C^\infty$-conjugated to Diophantine translations.
The proof of Theorem~\ref{thm:cf-vect-tori} is just a slight
modification of their proof.

\begin{proof}[Proof of Theorem~\ref{thm:cf-vect-tori}]
  Let $X(\theta)=(X_1(\theta),X_2(\theta),\ldots,X_d(\theta))$ be the
  coordinates of $X$ in the canonical trivialization of $T\bb{T}^d$
  and let $\Omega\in\Lambda^d(\bb{T}^d)$ be the only normalized
  $X$-invariant volume form given by
  Proposition~\ref{pro:invar-volume-form}. Let us define
  \begin{displaymath}
    \alpha_i\doteq\int_{\bb{T}^d}X_i\Omega\in\bb{R}, \quad\text{for }
    i=1,\ldots, d.
  \end{displaymath}
  So there exist smooth functions $u_i$ such that $\Lie_X u_i = - X_i
  + \alpha_i$. Then we can define a smooth map
  $f\colon\bb{T}^d\to\bb{T}^d$ writing
  \begin{displaymath}
    f(\theta)\doteq \theta + (u_1(\theta),u_2(\theta),\ldots,
    u_d(\theta)) \mod 1,\quad\forall\theta\in\bb{T}^d.  
  \end{displaymath}
  
  And then we have
  \begin{equation}
    \label{eq:Df-const}
    Df(X) = (X_i + \Lie_X u_i)_{i=1}^d =
    (\alpha_1,\alpha_2,\ldots,\alpha_d). 
  \end{equation}

  From equation~(\ref{eq:Df-const}) we can easily see that
  $f(\bb{T}^d)$ must be a coset of a closed connected subgroup of
  $\bb{T}^d$. By construction, $f$ is isotopic to the identity and so
  $f$ must be surjective. On the other hand, the set of critical
  points for $f$ is $\Phi_X$-invariant, and by Sard's theorem, it is
  not the whole torus.  Therefore, every point of $\bb{T}^d$ is
  regular and $f$ is a diffeomorphism, since tori do not admit any
  non-injective self-covering maps homotopic to the identity.

  As we already observed in Remark~\ref{rem:C-conj}, the set of
  cohomology-free vector fields is invariant by $C^\infty$-conjugacy.
  Hence, $X_\alpha=(\alpha_1,\alpha_2,\ldots,\alpha_d)$ must be
  cohomology-free too. Finally, it is rather easy to verify that then,
  $X_\alpha$ must be a Diophantine vector field (see \S 3.2.2
  in~\cite{katok-robinson}).
\end{proof}

\section{Topological Restrictions}
\label{sec:top-rest}

As it was already proved in Section~\ref{sec:stri-erg}, the flow
$\{\Phi_X^t\}$ induced by a cohomology-free vector field $X\in\X(M)$
is minimal and uniquely ergodic. In particular, $X$ cannot exhibit any
singularity, and so, the Euler characteristic of $M$ must vanish.

For a very long time this was the only known topological restriction
for manifolds supporting cohomology-free vector fields, until Federico
and Jana Rodr\'\i guez-Hertz produced a breakthrough in
\cite{fede-jana}, finding additional restrictions on the first Betti
number of the manifold.

For simplifying the exposition, let us first present a definition that
will be used all along this work:

\begin{definition}
  \label{def:good-fib}
  Given a closed $d$-manifold $M$ and a smooth vector field
  $Y\in\X(M)$, we say that a $p\colon M\to\bb{T}^n$ (where $n\leq d$)
  is a \emph{good fibration for } $Y$ if it is a smooth submersion and
  there exists a Diophantine vector $\alpha\in\bb{R}^n$ verifying
  $Dp(Y)\equiv\alpha$.
\end{definition}

Now we can state the main result of this section:

\begin{theorem}[F. \& J. Rodr\'\i guez-Hertz~\cite{fede-jana}]
  \label{thm:fede-jana}
  Let $X\in\X(M)$ be a cohomology-free vector field on the closed
  manifold $M$ and let us write $\beta_1\doteq\dim H_1(M,\bb{Q})$.
  Then there exists a good fibration $p\colon M\to\bb{T}^{\beta_1}$
  for $X$, where $Dp(X)\equiv X_\alpha\in\X(\bb{T}^{\beta_1})$ and
  $\alpha$ is a Diophantine vector. In particular, it holds
  $\beta_1(M)\leq\dim M$.
\end{theorem}

This fundamental result gives non-trivial information on the topology
of $M$ in all but one case: when $M$ has trivial first rational
homology group.

This is the main reason why it is necessary to attack Katok Conjecture
with different techniques, depending on the vanishing or not of the
first Betti number of the manifold.


\chapter{The case $\beta_1(M)\geq 1$}
\label{chap:b_1>0}

In this chapter we present the proof of Theorem A.

We continue assuming that $M$ is a closed orientable manifold and from
now on, we shall assume that $\dim M = 3$ and
\begin{displaymath}
  \beta_1(M)\doteq\dim H_1(M,\bb{Q})\geq 1. 
\end{displaymath}

For a better organization, we shall base the proof of Theorem A on the
following two propositions:

\begin{proposition}
  \label{pro:beta>1}
  Let us suppose that there exists $X\in\X(M)$ verifying:
  \begin{enumerate}
  \item The flow $\{\Phi_X^t\}_{t\in\bb{R}}$ induced by $X$ does not
    have any periodic orbit;
  \item and there is a good fibration $q\colon M\to\bb{T}^1$ for $X$.
  \end{enumerate}

  Then $\beta_1(M)\geq 2$.
\end{proposition}

\begin{proposition}
  \label{pro:beta=2}
  Let $X$ be a smooth vector field on $M$ and suppose that the induced
  flow $\{\Phi_X^t\}$ preserves a smooth volume form $\Omega$, \ie
  $\Lie_X\Omega\equiv 0$. Besides, assume that there exists a good
  fibration $p\colon M\to\bb{T}^2$ for $X$ verifying $Dp(X)=
  X_\alpha$.

  Then, if $M$ is not diffeomorphic to $\bb{T}^3$, $\Dis(X)$ has
  infinite dimension.
\end{proposition}

\section{Proof of Theorem A}
\label{sec:proof-thm-A}

This short section is devoted to prove Theorem A, assuming
Propositions~\ref{pro:beta>1} and~\ref{pro:beta=2}.

We are supposing that $M$ is a closed orientable $3$-manifold, with
$\beta_1(M)\geq 1$ and $X\in\X(M)$ is a cohomology-free vector field.
By Proposition~\ref{cor:minimality} we know that the induced flow
$\{\Phi_X^t\}$ is minimal, so in particular, it does not exhibit any
periodic orbit. On the other hand, by Theorem~\ref{thm:fede-jana}, we
know that there exists a good fibration $q\colon M\to\bb{T}^1$ for
$X$, with $Dq(X)$ verifying a Diophantine condition. Notice that in
the one-dimensional case, being Diophantine is equivalent to be
different from zero. Hence, we can apply Proposition~\ref{pro:beta>1}
for concluding that $\beta_1(M)\geq 2$.

Therefore, we can apply Theorem~\ref{thm:fede-jana} once again for
getting a good fibration $p\colon M\to\bb{T}^2$ for $X$ such that
$Dp(X)$ is a Diophantine vector in $\bb{R}^2$. On the other hand, by
Proposition~\ref{pro:invar-volume-form} we know that there exists a
smooth $X$-invariant volume form $\Omega$. And by
Proposition~\ref{pro:one-dim-invar-dist}, we can assure that $\dim
\Dis(X)=1$. So, if we apply Proposition~\ref{pro:beta=2}, we conclude
that $M$ is diffeomorphic to $\bb{T}^3$.

Finally, by Theorem~\ref{thm:cf-vect-tori}, $X$ is
$C^\infty$-conjugated to a constant vector field on $\bb{T}^3$, which
satisfies a Diophantine condition like
estimate~(\ref{eq:diophantine-def}), and we finish the proof of
Theorem A.

\section{Proof of Proposition~\ref{pro:beta>1}}
\label{sec:beta>1}

Let $X_\alpha\in\X(\bb{T}^1)$ be the Diophantine vector field given by
$X_\alpha \equiv Dq(X)$.  We know that $\alpha\neq 0$ and there is no
loss of generality supposing that $\alpha>0$.

Notice that for any $\theta\in\bb{T}^1$, the fiber $q^{-1}(\theta)$ is
a global transverse section for the flow $\{\Phi_X^t\}_{t\in\bb{R}}$.
So, it makes sense to define the Poincar\'e return map to
$q^{-1}(\theta)$ and this will be denoted by $\scr{P}_\theta$. Observe
that $\scr{P}_\theta = {\Phi_X}^{\alpha^{-1}}\Big|_{q^{-1}(\theta)}$.

Since the flow $\{\Phi_X^t\}$ does not have any periodic orbit, the
Poincar\'e return map $\scr{P}_\theta$ does not have any periodic
point. Hence, the Euler characteristic of the fiber $q^{-1}(\theta)$
must vanish.  Taking into account that the fiber is an orientable
(maybe non-connected) surface, we can affirm that it is diffeomorphic
to a disjoint union of $k$ $2$-torus. Our next step consists in
proving that we can modify our good fibration $q$ for getting another
one with connected fiber. This is the contents of our next

\begin{lemma}
  \label{lem:connec-fiber}
  If the fibration $q\colon M\to\bb{T}^1$ is such that
  $q^{-1}(\theta)$ has exactly $k$ connected components for some (and
  hence for any) $\theta\in\bb{T}^1$, then there exists another smooth
  good fibration $\tilde{q}\colon M\to\bb{T}^1$ satisfying:
  \begin{enumerate}
  \item $\tilde{q}^{-1}(\theta)$ is diffeomorphic to the $2$-torus;
  \item $D\tilde{q}(X)\equiv X_{k^{-1}\alpha}$;
  \item and the diagram
    \begin{equation*}
      \xymatrix{M\ar[rr]^{q}\ar[rd]_{\tilde{q}} && \bb{T}^1 \\
        &\bb{T}^1\ar[ru]_{E_k}&},
    \end{equation*}
    is commutative, where $E_k : \theta\mapsto k\theta$ is the
    canonical $k$-fold covering of the circle.
  \end{enumerate}
\end{lemma}

\begin{proof}
  Let $\theta_0$ be an arbitrary point of $\bb{T}^1$ and let us write
  $M_0$ for denoting a connected component of $q^{-1}(\theta_0)$.
  Since our manifold $M$ is connected, the Poincar\'e return map
  $\scr{P}_{\theta_0}$ must cyclically interchange all the connected
  components of $q^{-1}(\theta_0)$. Then, if we define
  \begin{equation*}
    M_t\doteq {\Phi_X}^{t\alpha^{-1}}(M_0),\quad\text{for every
    }t\in\bb{R},
  \end{equation*}
  it holds
  \begin{align*}
    M_t&=M_{t+k},\quad\text{for every }t\in\bb{R}; \\
    M&=\bigcup_{t\in\bb{R}}M_t.
  \end{align*}

  Therefore, if we define $\tilde{q}\colon M\to\bb{T}^1$ by
  \begin{equation*}
    \tilde{q}(x)\doteq k^{-1}t + \bb{Z}\in\bb{R}/\bb{Z},\quad\text{if }x\in M_t,
  \end{equation*}
  we easily see that $\tilde{q}$ is a good fibration for $X$, and it
  clearly satisfies properties (1), (2) and (3).
\end{proof}

With the purpose of simplifying our notation, we shall make the
assumption that our original good fibration $q\colon M\to\bb{T}^1$ was
such that its fibers $q^{-1}(\theta)$ were connected, and hence,
diffeomorphic to $\bb{T}^2$.

Now, let us fix a point $\theta_0\in\bb{T}^1$ and let $f\colon
q^{-1}(\theta_0)\to\bb{T}^2$ denote any diffeomorphism. Hence, we can
use the diffeomorphism $f$ to write the Poincar\'e return map
$\scr{P}_{\theta_0}$ as a diffeomorphism of $\bb{T}^2$, \ie we have
$f\circ\scr{P}_{\theta_0}\circ f^{-1}\in\Diff(\bb{T}^2)$. Then we can
choose an appropriate matrix $\hat{A}\in\SL(2,\bb{Z})$ such that its
induced linear automorphism $A\in\Diff(\bb{T}^2)$ is isotopic to
$f\circ\scr{P}_{\theta_0}\circ f^{-1}$.

By Lefschetz fixed point theorem, and since $\scr{P}_{\theta_0}$ is
fixed-point free, we know that
\begin{equation}
  \label{eq:lef-number}
  0=L(\scr{P}_{\theta_0})=\det(\hat{A}-id_{\bb{R}^2}).
\end{equation}

In this way, since $\hat{A}\in\SL(2,\bb{Z})$,
equation~(\ref{eq:lef-number}) implies that $1$ is the only element in
the spectrum of $\hat{A}$.  Therefore, $\hat{A}$ must be
$\SL(2,\bb{Z})$-conjugated to a matrix of the following form:
\begin{equation}
  \label{eq:A-jordan-form}
  \begin{pmatrix}
    1 & 0 \\
    n_0 & 1
  \end{pmatrix},
\end{equation}
commonly named the Jordan form of $\hat{A}$.

Then, post-composing $f$ with an appropriate element of
$\SL(2,\bb{Z})$ if necessary, we can assume that $\hat{A}$ equals to
matrix~(\ref{eq:A-jordan-form}).

On the other hand, notice that since $\scr{P}_{\theta_0}$ is the
time-$\alpha^{-1}$ map of the flow $\{\Phi_X^t\}$, matrix $\hat{A}$
(in fact, the conjugacy class of $\hat{A}$ in $\SL(2,\bb{Z})$)
determines the topology of $M$. More precisely, we know that $M$ is a
$\bb{T}^2$-bundle over $\bb{T}^1$, and so there exists a matrix
$\hat{B}\in\SL(2,\bb{Z})$ such that $M$ is smoothly diffeomorphic to
$\bb{T}^2\times\bb{R}\big/(B,1)$, where
$(B,1)\in\Diff(\bb{T}^2\times\bb{R})$ is defined by $(B,1)\colon
(x,t)\mapsto (B x,t-1)$. Furthermore, it is well-known that, given
$\hat{B}_1,\hat{B}_2\in\SL(2,\bb{Z})$,
$\bb{T}^2\times\bb{R}\big/(B_1,1)$ is homeomorphic to
$\bb{T}^2\times\bb{R}\big/(B_2,1)$ if and only if $\hat{B}_1$ and
$\hat{B}_2$ are $\SL(2,\bb{Z})$-conjugated (see for instance
\cite{hatcher}, Theorem 2.6, p. 36). Taking this into account, it is
not difficult to verify that $\hat{A}$, the only matrix which induces
an automorphism in the isotopy class of $f\circ\scr{P}_{\theta_0}\circ
f^{-1}$, and $\hat{B}$ must be conjugated in $\SL(2,\bb{Z})$.

Therefore, there exists a smooth diffeomorphism \mbox{$\Gamma\colon
  \bb{T}^2\times\bb{R}\big/(A,1)\to M$}.

Having gotten this nice topological characterization of $M$, our next
aim consists in studying the algebraic properties of the fundamental
group $\pi_1(M)$. For this, we define diffeomorphisms
$\tau_0,\tau_1,\tau_2\colon\bb{R}^3\to\bb{R}^3$ by
\begin{align}
  \tau_0 &: (x^0,x^1,x^2)\mapsto (x^0-1,x^1,x^2); \\
  \tau_1 &: (x^0,x^1,x^2)\mapsto (x^0,x^1-1,x^2); \\
  \tau_2 &: (x^0,x^1,x^2)\mapsto (x^0,x^1+n_0x^0,x^2-1).
\end{align}

If $G(\tau_i)$ denotes the subgroup of $\Diff(\bb{R}^3)$ generated by
$\{\tau_0,\tau_1,\tau_2\}$, we easily see that
\begin{displaymath}
  \bb{R}^3\Big/G(\tau_i) = \bb{T}^2\times\bb{R}\Big/(\hat{A},1),
\end{displaymath}
an therefore, we have that $G(\tau_i)$ is (algebraically) isomorphic
to $\pi_1(M)$. As a consequence of this, we have that $H_1(M,\bb{Q})$
is isomorphic to
\begin{displaymath}
  \Big(G(\tau_i)\Big/[G(\tau_i),G(\tau_i)]\Big)\otimes\bb{Q},
\end{displaymath}
where $[G(\tau_i),G(\tau_i)]$ denotes the commutator subgroup of
$G(\tau_i)$.

Hence, we finish the proof of Proposition~\ref{pro:beta>1} with the
following

\begin{lemma}
  It holds
  \begin{displaymath}
    \rank\Big(G(\tau_i)\Big/[G(\tau_i),G(\tau_i)]\Big)\geq 2.
  \end{displaymath}
\end{lemma}

\begin{proof}
  Let $H\doteq\mathrm{span}\{\tau_0,\tau_2\}$ be the subgroup of
  $G(\tau_i)$ generated by $\tau_0$ and $\tau_2$.  Let us write
  $\mathrm{pr}_i\colon\bb{R}^3\to\bb{R}$ for the canonical projection
  on the $i$-th coordinate, where $i=0,1,2$.

  First, notice that for any $g\in [G(\tau_i),G(\tau_i)]$, we have
  \begin{equation}
    \label{eq:even-disp-comm}
    \mathrm{pr}_i\circ g - \mathrm{pr}_i \equiv 0 , \quad
    \text{for } i=0,2.  
  \end{equation}

  Secondly, observe that $H$ is isomorphic to $\bb{Z}\oplus\bb{Z}$,
  being a possible group isomorphism defined by
  \begin{equation}
    \label{eq:H-isomorph}
    h \mapsto (\mathrm{pr}_0\circ h - \mathrm{pr}_0,
    \mathrm{pr}_2\circ h - \mathrm{pr}_2).
  \end{equation}

  Finally, taking into account~(\ref{eq:even-disp-comm})
  and~(\ref{eq:H-isomorph}), we easily conclude that the restriction
  to $H$ of the canonical projection of $G(\tau_i)$ on its
  abelianization is injective. In other words,
  $G(\tau_i)\big/[G(\tau_i),G(\tau_i)]$ contains a subgroup isomorphic
  to $\bb{Z}\oplus\bb{Z}$.
\end{proof}

\section{Proof of Proposition~\ref{pro:beta=2}}
\label{sec:beta=2}

This is the last section of the current chapter and it is devoted to
proving Proposition~\ref{pro:beta=2}.

By hypothesis, there exists a good fibration $p\colon M\to\bb{T}^2$
for $X$. Since $M$ is a closed $3$-manifold fibering over $\bb{T}^2$,
the fibers of $p$ must be diffeomorphic to the union of $k$ copies of
$\bb{T}^1$.  If $k>1$, then the idea is that we can apply
Lemma~\ref{lem:connec-fiber} ``twice'' to get a new good fibration
with connected fibers. This is what we are going to do first:

\begin{lemma}
  \label{lem:connec-fiber-2}
  There exists another good fibration $\tilde{p}\colon M\to\bb{T}^2$
  for $X$ verifying the following conditions:
  \begin{enumerate}
  \item $\tilde{p}^{-1}(\theta)$ is connected (and then, diffeomorphic
    to $\bb{T}^1)$, for every $\theta\in\bb{T}^2$.
  \item There exists $k_0,k_1\in\bb{N}$, such that
    $D\tilde{p}(X)\equiv X_{\tilde{\alpha}}$, where
    $\tilde{\alpha}\doteq(k_0^{-1}\alpha_0, k_1^{-1}\alpha_1)$.
  \end{enumerate}
\end{lemma}

\begin{remark}
  \label{rem:diop-alpha-tilde}
  A very simple but fundamental observation for the future is that the
  new vector $\tilde{\alpha}$ continues to be Diophantine. This can be
  simply proved observing that
  \begin{displaymath}
    \Big|r(k_0^{-1}\alpha_0)+s(k_1^{-1}\alpha_1)\Big| \geq
    \frac{C}{(\max\{\abs{rk_0^{-1}},\abs{sk_1^{-1}}\})^\tau} \geq
    \frac{C(\min\{k_0,k_1\})^\tau}{(\max\{\abs{r},\abs{s}\})^\tau}
  \end{displaymath}
\end{remark}

\begin{proof}[Proof of Lemma~\ref{lem:connec-fiber-2}]
  Heuristically, we could apply twice the method used in the proof of
  Lemma~\ref{lem:connec-fiber} to ``unfold'' the good fibration $p$
  along each direction of $\bb{T}^2$. Nevertheless, here we shall
  develop a different technique that makes the proof a little clearer.

  Let us start noticing that the fibration $p\colon M\to\bb{T}^2$
  induces a smooth foliation $\scr{F}$ on $M$ which leaves are the
  connected components of the fibers of $p$. Since $\scr{F}$ is a
  foliation with all its leaves compact, the space of leaves of
  $\scr{F}$, which will be denoted by $M/\scr{F}$, is a Hausdorff
  surface.

  Moreover, $p\colon M\to\bb{T}^2$ clearly factors through
  $M/\scr{F}$, \ie if $\tilde{p}_0\colon M\to M/\scr{F}$ denotes the
  canonical quotient map, then there exists a continuous map
  $p^\prime\colon M/\scr{F}\to\bb{T}^2$ making the following diagram
  commutative:
  \begin{equation*}
    \xymatrix{M\ar[rr]^{p}\ar[rd]_{\tilde{p}_0} && \bb{T}^2 \\
      &M/\scr{F}\ar[ru]_{p^\prime}&}
  \end{equation*}

  Then, we can easily see that $p^\prime\colon M/\scr{F}\to\bb{T}^2$
  is a $k$-fold covering map ($k$ is the number of connected
  components of any fiber of $p$) and therefore, $M/\scr{F}$ must be
  homeomorphic to $\bb{T}^2$. So, we can find two integers
  $k_0,k_1\in\bb{N}$, with $k=k_0k_1$, and a homeomorphism $h\colon
  M/\scr{F}\to\bb{T}^2$ expanding the previous diagram and getting
  \begin{equation}
    \label{eq:graph-M/F}
    \xymatrix{
      M\ar[rr]^{p}\ar[d]_{\tilde{p}_0} && \bb{T}^2 \\
      M/\scr{F}\ar[rr]^{h}\ar[urr]^{p^\prime} &&
      \bb{T}^2\ar[u]_{E_{k_0,k_1}}
    }
  \end{equation}
  where $E_{k_0,k_1} :
  (\theta^0,\theta^1)\mapsto(k_0\theta^0,k_1\theta^1)$ is a $k$-fold
  covering. In this way, the map $\tilde{p}\doteq h\circ\tilde{p}_0$
  is a smooth fibration satisfying
  $D\tilde{p}(X)\equiv(k_0^{-1}\alpha_0,k_1^{-1}\alpha_1)$ as desired.
\end{proof}

Having proved that we can find a good fibration with connected fibers
satisfying all the hypotheses of Proposition~\ref{pro:beta=2}, to
simplify the notation, we shall assume that our original good
fibration $p\colon M\to\bb{T}^2$ has connected fibers.

So, writing $q_0\doteq\mathrm{pr}_0\circ p\colon M\to\bb{T}^1$, we get
a good fibration for $X$ over $\bb{T}^1$, having connected fibers
diffeomorphic to $\bb{T}^2$. On the other hand, since $Dp(X)\equiv
X_\alpha$, being $\alpha$ a Diophantine vector of $\bb{R}^2$, we know
that $\{\Phi_X^t\}$ cannot exhibit any periodic orbit. Hence, we are
within the same context of Proposition~\ref{pro:beta>1}. Repeating the
same arguments exposed there, we may ensure that our manifold $M$ is
smoothly diffeomorphic to $\bb{T}^2\times\bb{R}\big/(A,1)$, where
\begin{equation}
  \label{eq:def-matrix-A}  
  \hat{A}\doteq
  \begin{pmatrix}
    1 & 0 \\
    n_0 & 1
  \end{pmatrix},
\end{equation}
and $(A,1)\in\Diff(\bb{T}^2\times\bb{R})$ is defined by $(A,1) :
(x,t)\mapsto (A x,t-1)$.

In this way we can reformulate the conclusion of
Proposition~\ref{pro:beta=2} saying that $\Dis(X)$ has infinite
dimension, provided that $n_0\neq 0$.

Continuing with the notation introduced in Section~\ref{sec:beta>1},
the Poincar\'e return map to the fiber $q_0^{-1}(\theta)$ shall be
denoted by $\scr{P}_\theta$, \ie $\scr{P}_\theta=
\Phi_X^{\alpha_0^{-1}}\Big|_{q_0^{-1}(\theta)}$.

Observe that all the things that we have done so far had the purpose
of returning to the setting of Proposition~\ref{pro:beta>1}.
Nevertheless, in this context we have additional geometric information
about the Poincar\'e return map $\scr{P}_\theta$. First, it preserves
a smooth foliation where all the leaves are circles (see
(\ref{eq:fol-F}) for the definition of the invariant foliation). As we
will see in paragraph~\ref{sec:inv-fol}, this will let us get a fine
system of coordinates for $\scr{P}_\theta$. Secondly, since
$\{\Phi_x^t\}$ preserves a smooth volume form, we easily see that
$\scr{P}_\theta$ also preserves a smooth volume form. We shall use
this in paragraph~\ref{sec:inv-vol-form} to improve our system of
coordinates proving that, in fact, $\scr{P}_\theta$ is linearizable,
\ie it is smoothly conjugated to an affine map in $\bb{T}^2$.

Finally, using the fact that $\scr{P}_\theta$ is $C^\infty$-conjugated
to an affine map and applying a classical construction, attributed to
Anatole Katok~\cite{katok-robinson}, we shall prove in
paragraph~\ref{sec:inv-distr} that there exist infinitely many linear
independent $X$-invariant distributions, provided $\scr{P}_\theta$ is
not isotopic to the identity.

\subsection{The Invariant Foliation}
\label{sec:inv-fol}

In this paragraph we shall use the fact that $\scr{P}_\theta\in
\Diff(q_0^{-1}(\theta))$ preserves a smooth foliation for proving that
$\scr{P}_\theta$ is smoothly conjugated to a skew-product over a rigid
rotation of $\bb{T}^1$.

For this, first notice that the flow $\{\Phi_X^t\}$ preserves the
codimension-two foliation in $M$ induced by the fibers of $p$.  In
fact it holds
\begin{equation}
  \label{eq:p-fibers-fol}
  \Phi_X^t(p^{-1}(\theta^0,\theta^1))=
  p^{-1}(\theta^0+t\alpha_0,\theta^1+t\alpha_1),
  \quad\forall(\theta^0,\theta^1)\in\bb{T}^2,\ \forall t\in\bb{R}. 
\end{equation}
Moreover, by definition, each fiber of $p$ is contained in a fiber of
$q_0$. In other words, the fibration $p$ is inducing a codimension-one
foliation on each fiber $q_0^{-1}(\theta)$, and this foliation happens
to be $\scr{P}_\theta$-invariant.

Then, let us fix some point $\theta\in\bb{T}^1$ and consider any
smooth diffeomorphism $f_0\colon\bb{T}^2\to q_0^{-1}(\theta)$. To
simplify forthcoming notation, let us define
\begin{displaymath}
  \scr{P}_1\doteq f_0^{-1}\circ\scr{P}_\theta\circ
  f_0\in\Diff(\bb{T}^2).
\end{displaymath}

Let $\mathcal{F}$ be the codimension-one foliation on $\bb{T}^2$
defined by
\begin{equation}
  \label{eq:fol-F} 
  \mathcal{F}(x)\doteq f_0^{-1}(p^{-1}(p(f_0(x)))), \quad\forall
  x\in\bb{T}^2, 
\end{equation}
where $\mathcal{F}(x)$ denotes the leaf of $\mathcal{F}$ passing
through $x$.

On the other hand, if we define the \emph{vertical foliation}
$\mathcal{V}$ in $\bb{T}^2$ by
\begin{equation}
  \label{eq:def-vert-fol}
  \mathcal{V}(\theta^0,\theta^1)\doteq \{\theta^0\}\times\bb{T}^1,
\end{equation}
and since all the leaves of $\mathcal{F}$ are diffeomorphic to
$\bb{T}^1$, it is a very well-known fact that there exists
$f_1\in\Diff(\bb{T}^2)$ verifying
\begin{equation}
  \label{eq:vert-fol}
  f_1(\mathcal{V}(x))=
  \mathcal{F}(f_1(x)), \quad\forall x\in\bb{T}^2.
\end{equation}

Once again, for the sake of simplicity, let us define
\begin{displaymath}
  \scr{P}_2\doteq f_1^{-1}\circ \scr{P}_1\circ f_1.
\end{displaymath}

From (\ref{eq:fol-F}) and (\ref{eq:vert-fol}) we easily see that there
exists $g_1\in\Diff(\bb{T}^1)$ satisfying
\begin{equation}
  \label{eq:lineariztion-x}
  \mathrm{pr}_0(\scr{P}_2(\theta^0,\theta^1))=g_1(\theta^0),
  \quad\forall(\theta^0,\theta^1)\in\bb{T}^2.
\end{equation}

Then we have the following

\begin{lemma}
  $g_1$ is smoothly linearizable, \ie there exists an
  orientation-preserving smooth diffeomorphism
  $h_1\colon\bb{T}^1\to\bb{T}^1$ such that $h_1^{-1}\circ g_1\circ
  h_1$ is an irrational rigid rotation on $\bb{T}^1$.
\end{lemma}

\begin{proof}
  First observe $g_1$ preserves orientation on $\bb{T}^1$, and hence
  it makes sense to consider its rotation number
  $\rho(g_1)\in\bb{T}^1$. Since $\scr{P}_2$ is a minimal
  diffeomorphism on $\bb{T}^2$, we have that $g_1$ does not exhibit
  any periodic point, and therefore, the rotation number $\rho(g_1)$
  is irrational and hence, it is completely determined by the order of
  the points of any orbit.

  Then, notice that the order of the points of
  $\{g_1^n(x)\}_{n\in\bb{Z}}$ in $\bb{T}^1$, for any $x$ in
  $\bb{T}^1$, is the same that the order of the leaves
  $\{\mathcal{F}(\scr{P}_\theta^n(z))\}_{n\in\bb{Z}}$ in $\bb{T}^2$,
  for any $z\in\bb{T}^2$.  On the other hand, we know the order of the
  leaves $\{\mathcal{F}(\scr{P}_\theta^n(z))\}_{n\in\bb{Z}}$ is given
  by the Poincar\'e return map to the global section
  $\{\theta_0\}\times\bb{T}^1\subset\bb{T}^2$ of the flow on
  $\bb{T}^2$ induced by the constant vector field
  $(\alpha_0,\alpha_1)$. We can easily see that the dynamics of this
  return map is given by the rigid rotation $x\mapsto
  x+\alpha_1/\alpha_0$.  Therefore, we can affirm that
  $\rho(g_1)=\alpha_1/\alpha_0 \mod \bb{Z}$.

  Besides, we know that, by hypothesis, there exist real positive
  constants $C$ and $\tau$ verifying
  \begin{displaymath}
    \abs{m\alpha_0+n\alpha_1}\geq
    \frac{C}{(\max\{\abs{m},\abs{n}\})^\tau},\quad \forall
    (m,n)\in\bb{Z}^2\setminus\{(0,0)\},
  \end{displaymath}
  and thus, elementary computations show that, indeed, it holds
  \begin{equation}
    \label{eq:diof-for-diff}
    \Big|m+n\frac{\alpha_1}{\alpha_0}\Big|\geq
    \frac{C^\prime}{\abs{n}^\tau}  \quad\forall
    n\in\bb{Z}\setminus\{0\}, 
  \end{equation}
  for some other real constant $C^\prime>0$.

  Finally, taking into account (\ref{eq:diof-for-diff}), we can apply
  Yoccoz linearization theorem \cite{yoccoz} to guarantee that $g_1$
  is smoothly conjugated to the rigid rotation
  $R_{\alpha_1/\alpha_0}$.
\end{proof}
 
This diffeomorphism $h_1$ can be used for defining
$f_2\in\Diff(\bb{T}^2)$ by $f_2 : (\theta^0,\theta^1)\mapsto
(h_1(\theta^0),\theta^1)$, getting as result
\begin{equation}
  \label{eq:first-linearization}
  f_2^{-1}(\scr{P}_2(f_2((\theta^0,\theta^1))=
  \left(\theta^0+\frac{\alpha_1}{\alpha_0},
    \theta^1 + n_0\theta^0 + \eta(\theta^0,\theta^1)\right), 
\end{equation}
for some $\eta\in C^\infty(\bb{T}^2,\bb{R})$ and for every
$(\theta^0,\theta^1)\in\bb{T}^2$.

Once again let us write
\begin{equation}
  \label{eq:P_3}
  \scr{P}_3\doteq f_2^{-1}\circ\scr{P}_2\circ f_2.
\end{equation}

\subsection{The invariant Volume Form}
\label{sec:inv-vol-form}

In this paragraph we show that there exists a smooth
$\scr{P}_\theta$-invariant volume form and analyze the consequences of
this.

By hypothesis we know that there exists a smooth $X$-invariant volume
form $\Omega\in\Lambda^3(M)$. So, if we write
\begin{equation}
  \label{eq:def-omega}
  \omega\doteq i_X\Omega,
\end{equation}
we get an $X$-invariant $2$-form. And since $X$ is transverse to $\ker
Dq_0$, we easily see that $\omega\big|_{q_0^{-1}(\theta)}$ is a
$\scr{P}_\theta$-invariant area form on $q_0^{-1}(\theta)$.

Therefore, defining
\begin{displaymath}
  \omega_3 \doteq (f_0\circ f_1\circ f_2)^*\omega\in\Lambda^2(\bb{T}),
\end{displaymath}
we get a $\scr{P}_3$-invariant area form. Making some abuse of
notation we can consider $\omega_3$ as an element of
$\M(\scr{P}_3)\subset\M(\bb{T}^2)$, identifying the area form with the
Borel finite measure that it induces on $\bb{T}^2$. Then, by
(\ref{eq:first-linearization}), we know that it holds
\begin{equation}
  \label{eq:proj-omega-3}
  (\mathrm{pr}_0)_*\omega_3=K\mathrm{Leb}^1,
\end{equation}
where $K\doteq\int_{\bb{T}^2}\omega_3$ is a positive real constant and
$\mathrm{Leb}^1$ denotes the Haar measure on $\bb{T}^1$.

At this point it would be desirable to know that the invariant measure
$\omega_3$ is a constant multiple of $\mathrm{Leb}^2$, the Haar
measure of $\bb{T}^2$. We could easily achieve our goal applying the
classical Moser's isotopy theorem \cite{moser}, but \emph{a priori} we
could not continue to have the skew-product structure of our
diffeomorphism. This is the reason why it is necessary to get a
``foliated version'' of Moser's isotopy theorem. The following can be
considered a two-dimensional reformulation of a more general result
due to Richard Luz and Nathan dos Santos~\cite{luz-santos}:

\begin{theorem}
  \label{thm:fol-moser}
  Let $\Omega_1,\Omega_2\in\Lambda^2(\bb{T}^2)$ be two volume forms
  and suppose they satisfy:
  \begin{equation*}
    \int_{\bb{T}^2}\Omega_1=\int_{\bb{T}^2}\Omega_2,\quad\text{and}
    \quad\Omega_1(\mathrm {pr}_0^{-1}(C))=\Omega_2(\mathrm{pr}_0^{-1}(C)),
  \end{equation*}
  for every Borel measurable set $C\subset\bb{T}^1$, where we are
  considering $\Omega_1$ and $\Omega_2$ as elements of $\M(\bb{T}^2)$.
  Then there exists $H\in\Diff(\bb{T}^2)$ isotopic to the identity
  verifying
  \begin{equation*}
    H^*\Omega_1=\Omega_2,\quad\text{and}\quad
    H(\mathcal{V}(x))=\mathcal{V}(H(x)), \quad\forall x\in\bb{T}^2,
  \end{equation*}
  where $\mathcal{V}$ is the vertical foliation in $\bb{T}^2$ defined
  in (\ref{eq:def-vert-fol}).
\end{theorem}

\begin{proof}
  See the proof of Theorem 6.1 in \cite{luz-santos}.
\end{proof}

Therefore, if we take into account (\ref{eq:proj-omega-3}),
Theorem~\ref{thm:fol-moser} lets us affirm that there exists a skew
product map $f_3\in\Diff(\bb{T}^2)$ verifying
\begin{displaymath}
  {f_3}^*\big(K(d\theta^0\wedge d\theta^1)\big)=\omega_3,
\end{displaymath}

From this we see that the diffeomorphism $\scr{P}_4\doteq
f_3^{-1}\circ\scr{P}_3\circ f_3\in\Diff(\bb{T}^2)$ preserves the Haar
measure and therefore, we can conclude that
\begin{displaymath}
  \scr{P}_4(\theta^0,\theta^1)=
  \left(\theta^0+\frac{\alpha_1}{\alpha_0}, 
    \theta^1+n_0\theta^0+\chi(\theta^0)\right),
\end{displaymath}
for some real function $\chi\in C^\infty(\bb{T}^1,\bb{R})$.

Since $\frac{\alpha_1}{\alpha_0}$ satisfies Diophantine
condition~(\ref{eq:diof-for-diff}), arguments analogous to those used
in Example~\ref{ex:diop-vect} let us prove that the rigid rotation
$x\mapsto x+\frac{\alpha_1}{\alpha_0}$ on $\bb{T}^1$ is
cohomology-free, and hence, we can find a function $\zeta\in
C^\infty(\bb{T}^1,\bb{R})$ verifying
\begin{displaymath}
  \zeta(x+\alpha_1\alpha_0^{-1}) - \zeta(x)=
  \chi(x)-\int_{\bb{T}^1}\chi\:\mathrm{d}(\mathrm{Leb}^1),\quad\forall
  x\in\bb{T}^1.   
\end{displaymath}

This function $\zeta$ can be used for linearizing $\scr{P}_4$. More
precisely, if we define $f_4 :(\theta^0,\theta^1)\mapsto
(\theta^0,\theta^1+\zeta(\theta^0))$, we get
\begin{equation}
  \label{eq:linear-P4}
  f_4^{-1}\Big(\scr{P}_4\big(f_4(\theta^0,\theta^1)\big)\Big)=
  \left(\theta^0+\frac{\alpha_1}{\alpha_0},
    \theta^1+n_0\theta^0+
    \int_{\bb{T}^1}\chi\:\mathrm{d}(\mathrm{Leb}^1)\right).  
\end{equation}

\subsection{Invariant Distributions}
\label{sec:inv-distr}

Summarizing what we have done in previous paragraphs, we can simply
say that there exists a diffeomorphism $F\colon\bb{T}^2\to
q_0^{-1}(\theta)$ verifying
\begin{equation}
  \label{eq:final-ver-linearized}
  F^{-1}\circ\scr{P}_\theta\circ F = A +
  (\alpha_1\alpha_0^{-1}, \beta),
\end{equation}
where $A$ is the automorphism of $\bb{T}^2$ induced by matrix
$\hat{A}$ defined in (\ref{eq:def-matrix-A}) and
$\beta=\int_{\bb{T}^1}\chi\:\mathrm{d}(\mathrm{Leb}^1)$ is obtained in
(\ref{eq:linear-P4}).

By (\ref{eq:def-matrix-A}), we know that if $n_0=0$, then $M$ is
diffeomorphic to $\bb{T}^3$. Hence, we shall assume that $n_0\neq 0$
and applying a construction due to Katok~\cite{katok-robinson}, we
will get infinitely many linearly independent
$\scr{P}_\theta$-invariant distributions on $\bb{T}^2$.

For this, let us start defining $T_m\in\Dis(\bb{T}^2)$, for each
$m\in\bb{Z}\setminus\{0\}$, writing
\begin{equation}
  \label{eq:def-inv-distrib}
  \langle T_m,\psi\rangle
  \doteq\sum_{k\in\bb{Z}}\hat{\psi}(kn_0m,m) e^{-2\pi
    ikm\left(\beta + \frac{k-1}{2}n_0\alpha_1\alpha_0^{-1}\right)},
\end{equation}
for each $\psi\in C^\infty(\bb{T}^2,\bb{R})$ and where
$\hat{\psi}\colon\bb{Z}^2\to\bb{C}$ denotes, as usual, the Fourier
transform of $\psi$. Clearly, the set
$\{T_m:m\in\bb{Z}\setminus\{0\}\}$ is linearly independent.
Furthermore, we can make the following

\begin{claim}
  If we define $B\doteq A+(\alpha_1\alpha_0^{-1},\beta)\in
  \Diff(\bb{T}^2)$, it holds
  \begin{equation}
    \label{eq:Tm-invariance}
    \langle T_m,\psi\circ B \rangle = \langle T_m, \psi\rangle,
    \quad\forall m\in\bb{Z}\setminus\{0\},\ \forall\psi\in
    C^\infty(\bb{T}^2,\bb{R}), 
  \end{equation}
\end{claim}

In fact, we have
\begin{equation*}
  \begin{split}
    \widehat{\psi\circ B}(k,\ell) & = \widehat{\psi\circ A}(k,\ell)
    \exp(2\pi i(k\alpha_1\alpha_0^{-1}+\ell\beta))  \\
    & = \hat{\psi}((A^*)^{-1}(k,\ell))\exp(2\pi
    i(k\alpha_1\alpha_0^{-1}+\ell\beta)) \\
    &= \hat{\psi}(k-n_0\ell,\ell)\exp(2\pi
    i(k\alpha_1\alpha_0^{-1}+\ell\beta)).
  \end{split}
\end{equation*}

And hence, it holds
\begin{equation*}
  \begin{split}
    \langle T_m ,\psi\circ B\rangle &= \sum_{k\in\bb{Z}}
    \widehat{\psi\circ B}(kn_0m,m) e^{-2\pi ikm
      \left(\beta+\frac{k-1}{2}n_0\alpha_1\alpha_0^{-1}\right)} \\
    &=\sum_{k\in\bb{Z}} \hat{\psi}((k-1)n_0m,m) e^{2\pi i
      (kn_0\alpha_1\alpha_0^{-1}+\beta)m} e^{-2\pi ikm
      \left(\beta+\frac{k-1}{2}n_0\alpha_1\alpha_0^{-1}\right)} \\
    &=\sum_{k\in\bb{Z}} \hat{\psi}((k-1)n_0m,m) e^{-2\pi i(k-1)m
      \left(\beta+\frac{k-2}{2}n_0\alpha_1\alpha_0^{-1}\right)} \\
    &=\langle T_m,\psi\rangle,
  \end{split}
\end{equation*}
for every $\psi\in C^\infty(\bb{T}^2,\bb{R})$.

Then, taking into account equation~(\ref{eq:final-ver-linearized}), we
see that we may push-forward each $T_m$ by $F$ for getting infinitely
many linearly independent $\scr{P}_\theta$-invariant distributions on
$q_0^{-1}(\theta)$.

Finally, if we define $\tilde{T}_m\in\Dis(M)$, for
$m\in\bb{Z}\setminus\{0\}$, writing
\begin{equation}
  \label{eq:def-inv-dist-M}
  \langle\tilde{T}_m,\psi\rangle\doteq \int_{\bb{T}^1} \Big\langle
  F_*T_m,(\psi\circ\Phi_X^{-t})\big|_{q_0^{-1}(\theta+t)}
  \Big\rangle\:\mathrm{d}t,
\end{equation}
for every $\psi\in C^\infty(M,\bb{R})$, we easily see that each
$\tilde{T}_m\in\Dis(X)\setminus\M(X)$ and they clearly form a linearly
independent set.


\chapter{The case $\beta_1(M)=0$}
\label{chap:b_1=0}

This chapter aims to prove that there is no cohomology-free vector
field on closed orientable $3$-manifolds with vanishing first Betti
number.

First of all, notice that by Poincar\'e duality, a closed $3$-manifold
with trivial first rational cohomology group must also have trivial
second cohomology group. So, from now on and until the end of the
current chapter, $M$ will denote a rational homological $3$-sphere and
we shall suppose that there exists a cohomology-free $X\in\X(M)$.

The general strategy for getting a contradiction from our assumptions
consists in proving first that there exists an $X$-invariant one-form
with no singularity. Then, we shall analyze the integrability of its
kernel, getting two possible cases: either the kernel of the invariant
form is everywhere integrable, or it is a contact structure, being $X$
collinear with the induced Reeb vector field (see
Section~\ref{sec:contact-struct-case} for more details). The rest of
the proof consists in proving that both cases lead to a contradiction.

\section{The Invariant $1$-form}
\label{sec:inv-one-form}

The main purpose of this section is to prove that the derivative of
the flow $\{\Phi_X^t\}$ preserves a smooth two-dimensional plane
field.

At this point the author would like to thank Giovanni Forni who kindly
communicated the following result to us:

\begin{theorem}
  \label{thm:inv-one-form}
  Let $M$ be a closed $3$-manifold such that
  $H^1(M,\bb{Q})=H^2(M,\bb{Q})=0$ and let $X\in\X(M)$ be a
  cohomology-free vector field. Then there exists
  $\lambda\in\Lambda^1(M)$ verifying
  \begin{equation*}
    \Lie_X\lambda\equiv 0 \quad\text{and}\quad \lambda(p)\neq 0, 
  \end{equation*}
  for every $p\in M$.
\end{theorem}

\begin{proof}
  By Proposition~\ref{pro:invar-volume-form}, we know that there
  exists an $X$-invariant volume form $\Omega\in\Lambda^3(M)$. Hence,
  if we write $\omega\doteq i_X\Omega$, Cartan's formula lets us
  affirm
  \begin{displaymath}
    0=\Lie_X\Omega = d(i_X\Omega)+i_X(d\Omega)=d\omega,
  \end{displaymath}
  \ie $\omega\in\Lambda^2(M)$ is a closed form. On the other hand, by
  the Universal Coefficient Theorem we know that $H^2(M,\bb{R})=0$,
  and thus, there exists a $1$-form $\tilde{\lambda}$ such that
  $\omega=d\tilde{\lambda}$.  Applying Cartan's formula once again we
  obtain
  \begin{displaymath}
    \Lie_X\tilde{\lambda}=d(i_X\tilde{\lambda})+i_X(d\tilde{\lambda})
    = d(i_X\tilde{\lambda})+i_X(i_X\Omega)=d(i_X\tilde{\lambda}). 
  \end{displaymath}

  Notice that $i_X\tilde{\lambda}$ is an element of
  $C^\infty(M,\bb{R})$, so there exists a smooth function $u\colon
  M\to\bb{R}$ verifying
  \begin{equation}
    \label{eq:cohomo-eq-lambda}
    \Lie_X u = - i_X\tilde{\lambda} + \int_M
    (i_X\tilde{\lambda})\Omega. 
  \end{equation}

  Therefore, if we define $\lambda\doteq \tilde{\lambda}+du$, it still
  holds $d\lambda=\omega$ and besides,
  \begin{equation*}
    \begin{split}
      \Lie_X\lambda & =\Lie_X\tilde{\lambda} + \Lie_X du =
      d(i_X\tilde{\lambda}) + d(i_X du) \\
      & = d\big(i_X\tilde{\lambda} + \Lie_X u \big) = d\bigg(\int_M
      (i_X\tilde{\lambda})\Omega \bigg) = 0,
    \end{split}
  \end{equation*}
  \ie $\lambda$ is an $X$-invariant $1$-form.

  Then, taking into account the minimality of $\{\Phi_X^t\}$, we
  easily see that $\lambda$ exhibits a singularity if and only if
  $\lambda\equiv 0$. On the other hand, since $d\lambda=i_X\Omega\neq
  0$, we know that $\lambda\not\equiv 0$, and therefore, $\lambda$
  does not have any singularity.
\end{proof}

So, we have proved the existence of a singularity-free $1$-form
$\lambda$ on $M$ which is invariant under the flow $\{\Phi_X^t\}$.
This lets us define an invariant two-dimensional plane field
\begin{displaymath}
  \Sigma\doteq \ker\lambda\subset TM.
\end{displaymath}

Now, it seems to be natural to ask ourselves about the integrability
of the plane field $\Sigma$. For this, it is interesting to notice
that the minimality of $\{\Phi_X^t\}$ implies that $\Sigma$ is either
a contact structure or it is everywhere integrable.

These cases will be analyzed in the following two sections.

\section{The Contact Structure Case}
\label{sec:contact-struct-case}

Let us start this section recalling some fundamental facts about
Contact Geometry.

Given a $(2n+1)$-manifold $N$, we say that $\alpha\in\Lambda^1(N)$ is
a \emph{contact form} if $\alpha\wedge(d\alpha)^{\wedge 2n}$ is a
volume form on $N$. This clearly implies that $\ker\alpha\oplus\ker
d\alpha=TM$, and consequently, there exists a unique vector field
$Y\in\X(N)$, called the \emph{Reeb vector field induced by $\alpha$},
verifying
\begin{displaymath}
  i_Y\alpha\equiv 1,\quad\text{and}\quad i_Yd\alpha\equiv 0.
\end{displaymath}

As the completely opposite case we know by Froebenius theorem that the
kernel of a singularity-free $1$-form $\beta$ is completely integrable
(\ie there exists a smooth foliation $\scr{F}$ verifying
$T\scr{F}=\ker\beta$) if and only if $\beta\wedge d\beta\equiv 0$.

As it was already mentioned at the end of Section
\ref{sec:inv-one-form}, we have a strict dichotomy: either $\lambda$
is a contact form, or $\Sigma$ is completely integrable. In this
section we shall analyze the first case, \ie we shall assume that
$\lambda$ is a contact form.

We know that $d\lambda=\omega=i_X\Omega$, and therefore, $\ker
d\lambda=\bb{R}X$. This implies that $X\not\in\Sigma$ and
consequently, by equation (\ref{eq:cohomo-eq-lambda}), we have
\begin{displaymath}
  \lambda(X)\equiv\int_M (i_X\tilde{\lambda})\Omega\neq 0.
\end{displaymath}

All this implies that $X$ is a constant multiple of the Reeb vector
field of $\lambda$, and in particular, they have the same orbits.

A very important problem in Contact Geometry that has received a lot
of attention and has led much of the research in this area during the
last decades is the following conjecture proposed by Alan Weinstein in
\cite{weinstein}:

\begin{conjecture}[Weinstein's Conjecture]
  Let $N$ be a closed $3$-manifold, $\alpha\in\Lambda^1(N)$ be a
  smooth contact form and $Y\in\X(N)$ be its Reeb vector field. Then
  $Y$ exhibits a periodic orbit.
\end{conjecture}

Clifford Taubes has recently proved the validity of this conjecture in
\cite{taubes}, and in our setting, this leads us to a contradiction:
since $\{\Phi_X^t\}$ is minimal, it cannot have any periodic orbit.

\section{The Completely Integrable Case}
\label{sec:comp-integrable-case}

In this section we shall analyze the situation where $\Sigma$ is a
completely integrable plane field. As it was already mentioned in
Section~\ref{sec:out-main-results}, while this work was in progress,
Giovanni Forni communicated to the author that he had been able to
exclude this case using the foliation tangent to $\Sigma$ to prove
that $M$ should be diffeomorphic to a nilmanifold and $\{\Phi_X^t\}$
smoothly conjugated to a homogeneous flow. On the other hand, Stephen
Greenfield and Nolan Wallach had already proved in \cite{hypo-vect}
that $\bb{T}^3$ was the only $3$-dimensional nilmanifold that supported
cohomology-free homogeneous vector fields\footnote{In fact, in
  \cite{hypo-vect} they proved this for globally hypoelliptic vector
  fields.} (see \cite{fla-for} for higher dimensional nilmanifolds).

Nevertheless, in this work we propose a completely different proof
that does not use the integrability condition in a direct form. In
fact, the only information we need for our proof is that our vector
field is contained in the plane field $\Sigma$. This approach has the
advantage that seems to be more ``usable'' for solving the contact
structure case independently of Taubes' proof of Weinstein's
conjecture, which would be very desirable (see
Section~\ref{sec:remarks-beta1-0} for a more detailed discussion about
this point).

Our general strategy consists in proving that, under our assumptions
about the topology of $M$, $\{\Phi_X^t\}$ must be a positively
expansive flow (see Definition~\ref{def:expansive-flow}).

For getting this, we will have to carefully study the dynamics of the
derivative of the flow $\{\Phi_X^t\}$ on $TM$. This analysis starts in
paragraph~\ref{sec:proj-flow}, where we get our first result about the
angular behavior of $D\Phi_X\colon TM\times\bb{R}\to TM$, studying the
dynamics of the \emph{projective flow} (see
paragraph~\ref{sec:normal-proj-flow} for definitions). Then, in
paragraph~\ref{sec:normal-flow}, we shall get some results about the
radial behavior of flow $\{D\Phi_X^t\}$ analyzing the dynamics of the
\emph{normal flow} and proving that it exhibits a parabolic behavior.
And in paragraph~\ref{sec:expansivness}, we will prove that our flow
$\{\Phi_X^t\}$ is indeed positively expansive.

On the other hand, using a nice result due to Miguel
Paternain~\cite{patern} about expansive flows on $3$-manifolds, we
shall prove in paragraph~\ref{sec:expansivness} that there is no
closed $3$-manifold supporting positively expansive flows, getting our
desired contradiction.

\subsection{The Normal and Projective Flows}
\label{sec:normal-proj-flow}

This short paragraph is devoted to introduce some terminology that we
shall repeatedly use in subsequent paragraphs.

Let us start defining the relation $\sim$ on $TM$ by
\begin{displaymath}
  v\sim w \iff \pi(v)=\pi(w)\text{ and } v-w\in\mathrm{span}(X),
\end{displaymath}
where $\pi\colon TM\to M$ stands for the canonical vector bundle
projection.  This is clearly an equivalence relation and thus, we can
define $NX$ to be the quotient of $TM$ by this relation. This set $NX$
can be naturally endowed with a unique $C^\infty$ vector bundle
structure $\pi_N\colon NX\to M$ such that the quotient map
$\mathrm{pr}_X\colon TM\to NX$ given by $\mathrm{pr}_X :
v\mapsto\hat{v}\doteq\{w\in TM : v\sim w \}$ is a smooth vector bundle
map. This shall be called the \emph{normal vector bundle induced by}
$X$.

Observe that since $D\Phi_X^t(X(p))=X(\Phi_X^t(p))$, we easily see the
derivative of $\{\Phi_X^t\}$ induces a \emph{vector bundle flow}
$N\Phi_X\colon NX\times\bb{R}\to NX$ over $\{\Phi_X^t\}$, \ie it makes
sense to define
\begin{displaymath}
  N\Phi_X^t(\hat{v})\doteq \mathrm{pr}_X(D\Phi_X^t(v)),\quad\text{for
    any } v\in\mathrm{pr}_X^{-1}(\hat{v}),
\end{displaymath}
and any $t\in\bb{R}$. This flow $\{N\Phi_X^t\}$, which will be called
the \emph{normal flow} induced by $\{\Phi_X^t\}$, clearly verifies
$\pi_N\circ N\Phi_X^t=\Phi_X^t\circ\pi_N$, being $N\Phi_X^t\colon
NX_p\to NX_{\Phi_X^t(p)}$ a linear isomorphism.

Then, since $\{N\Phi_X^t\}$ is a vector bundle flow, it induces a new
flow on $\pi_\bb{P}\colon \bb{P}(NX)\to M$, the projectivization of
the normal bundle $\pi_N\colon NX\to M$. This will be called the
\emph{projective flow} induced by $\{\Phi_X^t\}$ and it shall be
denoted by $P\Phi_X\colon \bb{P}(NX)\times\bb{R}\to\bb{P}(NX)$. We
will write $\mathrm{pr}_{\bb{P}}\colon
NX\setminus\{0\}\to\bb{P}(NX)$\footnote{In this context $\{0\}$ means
  the zero section of $NX$.}  for the canonical quotient projection
given by $\mathrm{pr}_{\bb{P}} : \hat{v}\mapsto
(\bb{R}\setminus\{0\})\hat{v}$.

\subsection{Dynamics of the Normal Flow I}
\label{sec:normal-flow}

In this paragraph we start the analysis of the dynamics of the normal
flow $\{N\Phi_X^t\}$.

Let us start introducing any smooth Riemannian structure
$\langle\cdot,\cdot\rangle$ on $TM$. This naturally induces another
Riemannian structure $\langle\cdot,\cdot\rangle_{NX}$ on $NX$
defining, for each $p\in M$,
\begin{displaymath}
  \langle v_1,v_2\rangle_{NX}\doteq \langle
  v^\prime_1,v^\prime_2\rangle, \quad\forall v_1,v_2\in NX_p,
\end{displaymath}
where $v_i^\prime$ is defined as the only element of $T_pM$ verifying
simultaneously $\langle X(p),v_i^\prime(p)\rangle=0$ and
$\mathrm{pr}_X(v_i^\prime)=v_i$. The Finsler structures induced by
$\langle\cdot,\cdot\rangle$ and $\langle\cdot,\cdot\rangle_{NX}$ will
be denoted by $\|\cdot\|$ and $\|\cdot\|_{NX}$, respectively. As
usual, we shall also use the Riemannian structures
$\langle\cdot,\cdot\rangle$ and $\langle\cdot,\cdot\rangle_{NX}$ for
measuring angles between non-null vectors of the same fiber. Making
some abuse of notation, we shall use the symbol
$\sphericalangle(\cdot,\cdot)$ for both.

Before we state our first result about the radial behavior of vectors
in $NX$, we need to recall some notions of hyperbolic dynamics:

Given a closed $d$-manifold $B$, a vector bundle $\pi\colon E\to B$
and a non-singular vector field $Y\in\X^r(B)$ ($r\geq 2$), we say that
$A\colon E\times\bb{R}\to E$ is a \emph{linear cocycle over}
$\{\Phi_Y^t\}$ if it holds $\Phi_X^t\circ\pi=\pi\circ A(\cdot,t)$, for
any $t$, being the maps
$A(\cdot,t)\colon\pi^{-1}(p)\to\pi^{-1}(\Phi_X^t(p))$ linear
isomorphisms that verify
\begin{displaymath}
  A(p,t_0+t_1)=A(\Phi_X^{t_0}(p),t_1)A(p,t_0),\quad\forall p\in B,\
  \forall t_0,t_1\in\bb{R}.
\end{displaymath}

We shall say that cocycle $A$ is \emph{Anosov} if there exist two
sub-bundles $E^s,E^u\subset E$ and real constants $C>0$ and
$\rho\in(0,1)$ verifying
\begin{itemize}
\item $E^s\oplus E^u=E$,
\item $A(E^\sigma_p,t)=E^\sigma_{\Phi_Y^t(p)}$, for every $p\in M$,
  every $t\in\bb{R}$ and $\sigma=s,u$.
\item $\left\|A(\cdot,t)\big|_{E^s}\right\|\leq C\rho^t$, and
  $\left\|A(\cdot,-t)\big|_{E^u}\right\|\leq C\rho^t$, for any $t>0$,
\end{itemize}
where $\|\cdot\|$ is any Finsler structure on $\pi\colon E\to B$.

We shall say that cocycle $A$ is \emph{quasi Anosov} if, given any
$v\in E$, it holds
\begin{equation}
  \label{eq:def-quasi-Anosov}
  \sup_{t\in\bb{R}} \left\| A(v,t) \right\| <\infty \Rightarrow v=0. 
\end{equation}

The following result appears in different forms, and in fact with
different hypothesis, in the works of Ricardo Ma\~n\'e~\cite{mane},
Robert Sacker and George Sell~\cite{sacker}, and James
Selgrade~\cite{selgrade,selgrade-err}:

\begin{proposition}
  \label{pro:quasi-Anosov-impl-Anosov}
  If the flow $\{\Phi_Y^t\}$ does not have wandering points, then a
  cocycle $A$ is quasi Anosov if and only if it is Anosov.
\end{proposition}

On the other hand, we shall say that a vector field $Y\in\X^r(B)$ is
\emph{Anosov} if there exists a codimension-one $D\Phi_Y$-invariant
sub-bundle $F\subset TB$ verifying $F\oplus\bb{R}X=TB$ and such that
$D\Phi_Y|_F\colon F\times\bb{R}\to F$ is an Anosov linear cocycle.

Then, we can state the following result due to Claus Doering:

\begin{proposition}[Doering~\cite{doering}]
  \label{pro:Anosov-normal-impl-Anosov}
  Let suppose that $\{\Phi_Y^t\}$ does not have any wandering point.
  Then, $Y$ is an Anosov vector field if and only if its normal flow
  $\{N\Phi_Y^t\}$ is an Anosov linear cocycle (over $\{\Phi_Y^t\}$).
\end{proposition}

Now, we can present our first result about the dynamics of our normal
flow $\{N\Phi_X^t\}$:

\begin{lemma}
  \label{lem:quasi-anosov}
  There exists $\hat{v}_0\in NX$ such that $\hat{v}_0\neq 0$ and
  \begin{equation}
    \label{eq:quasi-anosov-condition}
    \sup_{t\in\bb{R}} \| N\Phi_X^t(\hat{v}_0)\|_{NX} < \infty.
  \end{equation}
\end{lemma}

\begin{proof}
  Let us suppose that estimate~(\ref{eq:quasi-anosov-condition}) is
  not satisfied by any non-vanishing vector in $NX$. In other words,
  let suppose that $N\Phi_X\colon NX\to\bb{R}\to NX$ is quasi Anosov.
  By Proposition~\ref{pro:quasi-Anosov-impl-Anosov}, $\{N\Phi_X^t\}$
  is an Anosov cocycle. Then,
  Proposition~\ref{pro:Anosov-normal-impl-Anosov} lets us affirm that
  $X$ is indeed Anosov.

  Finally, it is a very well-known fact that any Anosov flow exhibits
  (infinitely many) periodic orbits, which clearly contradicts the
  minimality of $\{\Phi_X^t\}$.
\end{proof}

Our second result about the dynamics of the normal flow is the
following

\begin{lemma}
  \label{lem:conserv-normal-flow}
  The normal flow $\{N\Phi_X^t\}$ is conservative. More precisely,
  there exists a symplectic form $\kappa$ on the vector bundle
  $\pi_N\colon NX\to M$ which is invariant under the action of
  $\{N\Phi_X^t\}$.
\end{lemma}

\begin{proof}
  Notice that $\omega=i_X\Omega=d\lambda$ is a $2$-form on $TM$
  verifying $i_X \omega\equiv 0$.  This implies that we may
  push-forward this form by $\mathrm{pr}_X$ on $NX$, \ie we can find a
  smooth $2$-form $\kappa$ on $NX$ such that
  \begin{displaymath}
    \kappa(\mathrm{pr}_X(v),\mathrm{pr}_X(w))=\omega(v,w), \quad\forall
    v,w\in T_p M,\ \forall p\in M.
  \end{displaymath}
  It is very easy to verify that $\kappa$ is symplectic on $NX$ and
  that it is $N\Phi_X$-invariant.
\end{proof}

\subsection{Dynamics of the Projective Flow}
\label{sec:proj-flow}

This paragraph aims to prove that the dynamics of the projective flow
is very simple. In fact, we shall get that the limit set of
$\{P\Phi_X^t\}$ is a smooth submanifold of $\bb{P}(NX)$ which happens
to be a graph over $M$, being the dynamics on this set smoothly
conjugated to $\{\Phi_X^t\}$.

For this, first we will need the following result due to Hiromichi
Nakayama and Takeo Noda about the geometry and amount of minimal sets
for the projective flow:

\begin{theorem}[Nakayama \& Noda~\cite{nakayama}]
  \label{thm:nakayama-inv-set}
  Let $V$ be a closed $3$-manifold and let $Y\in\X(V)$ be such that
  its induced flow $\Phi_Y\colon V\times\bb{R}\to V$ is minimal.

  Let $P\Phi_Y\colon \bb{P}(NY)\times\bb{R}\to\bb{P}(NY)$ be the
  projective flow induced by $\{\Phi_Y^t\}$. Hence, we have:
  \begin{enumerate}
  \item If $\{P\Phi_Y^t\}$ exhibits more than two minimal sets, then
    $V$ is diffeomorphic to $\bb{T}^3$ and $\{\Phi_X^t\}$ is
    continuously conjugate to an irrational translation.
  \item If $\{P\Phi_Y^t\}$ exhibits exactly two minimal sets
    $M_1,M_2\subset\bb{P}(NY)$ and $\{\Phi_X^t\}$ is not
    $C^0$-conjugate to an irrational translation on $\bb{T}^3$, then
    for any $z\in V$ it holds: $M_1\cap\pi_{\bb{P}}^{-1}(z)$ or
    $M_2\cap\pi_{\bb{P}}^{-1}(z)$ consists of a single point.
    Moreover, there exists a residual subset $B\subset V$ such that
    both sets $M_1\cap\pi_{\bb{P}}^{-1}(z)$ and
    $M_2\cap\pi_{\bb{P}}^{-1}(z)$ contain just a point, for every
    $z\in B$.
  \end{enumerate}
\end{theorem}

Since we are assuming that $H_1(M,\bb{Q})=0$,
Theorem~\ref{thm:nakayama-inv-set} lets us affirm that the flow
$\{P\Phi_X^t\}$ exhibits at most two minimal sets.

One is given by the plane field $\Sigma$. In fact, we have
$X(p)\in\Sigma_p$ for each $p\in M$, and hence,
\begin{equation}
  \label{eq:def-E-sigma}
  E_\Sigma\doteq\mathrm{pr}_X(\Sigma)\subset NX,
\end{equation}
is a smooth one-dimensional vector sub-bundle of $NX$. In this way,
$E_\Sigma$ determines exactly one point on each fiber of
$\pi_{\bb{P}}\colon\bb{P}(NX)\to M$. More precisely, we may define the
point $\theta_p\in\pi_{\bb{P}}^{-1}(p)$ by $\theta_p\doteq
\mathrm{pr}_{\bb{P}}({E_\Sigma}_p\setminus\{0\})$.

Notice that since the plane field $\Sigma$ is invariant under the
action of $\{D\Phi_X^t\}$, we have the flow $\{N\Phi_X^t\}$ leaves
invariant the line field $E_\Sigma$, and therefore, it holds
$P\Phi_X^t(\theta_p)=\theta_{\Phi_X^t(p)}$, for any $p\in M$ and any
$t\in\bb{R}$. So, summarizing we have that
\begin{equation}
  \label{eq:def-K-sigma}
  K_\Sigma\doteq\{\theta_p:p\in M\}\subset\bb{P}(NX)
\end{equation}
is a minimal set for $\{P\Phi_X^t\}$.

Finally, as it was mentioned above, we shall prove that $K_\Sigma$ is
indeed the only minimal set, and consequently, it is the $\alpha$- and
$\omega$-limit of any point in $\bb{P}(NX)$:

\begin{theorem}
  \label{thm:minimal-set-P(NX)}
  $K_\Sigma\subset\bb{P}(NX)$ defined in (\ref{eq:def-K-sigma}) is the
  only minimal set for $\{P\Phi_X^t\}$.
\end{theorem}

For proving Theorem~\ref{thm:minimal-set-P(NX)}, we shall suppose that
there exists another $P\Phi_X$-invariant minimal set
$K_0\subset\bb{P}(NX)$ (\ie different from $K_\Sigma$), and for the
sake of clarity of the exposition, we will separate the proof in
several lemmas:

\begin{lemma}
  \label{lem:orient-E-Sigma}
  Sub-bundle $E_\Sigma\subset NX$ defined in (\ref{eq:def-E-sigma}) is
  orientable, and therefore, it admits a non-vanishing section
  $\hat{Y}_0\in\Gamma(E_\Sigma)$.
\end{lemma}

\begin{proof}
  Since $\Sigma$ was defined as the kernel of a non-singular $1$-form
  and by hypothesis, $M$ is orientable, we have that $\Sigma\to M$ is
  orientable. On the other hand, our vector field $X$ can be
  considered as a non-singular element of $X\in\Gamma(\Sigma)$.

  This lets us affirm that $\Sigma\to M$ is a globally trivial vector
  bundle, and therefore, we can find a smooth section
  $Y_0\in\Gamma(\Sigma)$ verifying
  $\Sigma_p=\mathrm{span}\{X(p),Y_0(p)\}$, for every $p\in M$.

  Finally, defining $\hat{Y}_0\doteq\mathrm{pr}_X(Y_0)$ we get our
  desired section of $E_\Sigma\to M$.
\end{proof}

\begin{lemma}
  \label{lem:dynamics-on-E-Sigma}
  Assuming that there exists another minimal set
  $K_0\subset\bb{P}(NX)$, we can find a non-vanishing
  $\hat{Y}\in\Gamma(E_\Sigma)$ verifying
  \begin{equation}
    \label{eq:invariance-hat-Y}
    N\Phi_X^t\big(\hat{Y}(p)\big)=\hat{Y}(\Phi_X^t(p)),
    \quad\forall p\in M,\ \forall t\in\bb{R}. 
  \end{equation}
\end{lemma}

\begin{proof}
  Let $L_\Sigma\in C^\infty(M,\bb{R})$ be defined by
  \begin{displaymath}
    L_\Sigma(p)\hat{Y}_0(p)= \lim_{t\to 0}
    \frac{N\Phi_X^{-t}(\hat{Y}_0(\Phi_X^t(p)))-\hat{Y}_0(p)}{t},
    \quad\forall p\in M.
  \end{displaymath}
  
  Using the fact that $X$ is cohomology-free, we get a function $u\in
  C^\infty(M,\bb{R})$ verifying
  \begin{equation}
    \label{eq:lyapunov-exp-E-sigma-cohomo-eq}
    \Lie_X u = - L_\Sigma + \int_M L_\Sigma\Omega. 
  \end{equation}  

  Then, if we define $\hat{Y}\doteq e^u\hat{Y}_0$, applying
  equation~(\ref{eq:lyapunov-exp-E-sigma-cohomo-eq}) we clearly get
  \begin{displaymath}
    \lim_{t\to 0}
    \frac{N\Phi_X^{-t}(\hat{Y}(\Phi_X^t(p)))-\hat{Y}(p)}{t} =
    \left(\int_M L_\Sigma\Omega\right)\hat{Y}(p),
    \quad\forall p\in M,
  \end{displaymath}
  and therefore, it holds
  \begin{equation}
    \label{eq:Y-hat-growth}
    N\Phi_X^t(\hat{Y}(p))=\exp\bigg(t\int_M L_\Sigma\Omega\bigg)
    \hat{Y}(\Phi_X^t(p)), 
  \end{equation}
  for every $p\in M$ and every $t\in\bb{R}$.

  Notice that by equation~(\ref{eq:Y-hat-growth}), $\int_M
  L_\Sigma\Omega$ is a Lyapunov exponent of the linear cocycle
  $\{N\Phi_X^t\}$. So, let us suppose that $\int_M L_\Sigma\Omega\neq
  0$. In this case, the one-dimensional sub-bundle $E_\Sigma\subset
  NX$ is uniformly hyperbolic.

  On the other hand, by Theorem~\ref{thm:nakayama-inv-set}, we know
  that $K_0$ and $K_\Sigma$ are the only minimal sets on $\bb{P}(NX)$,
  and moreover, we can find a point $p_0\in M$ such that
  $\theta^\prime\in\bb{P}(NX)$ is the only point in
  $K_0\cap\pi_{\bb{P}}^{-1}(p_0)$.

  Observe that, since $K_0$ and $K_\Sigma$ are disjoint closed sets,
  we have that there exists a real constant $C>0$ such that
  \begin{equation}
    \label{eq:theta_p-theta-0-dist}
    \mathrm{dist}_{\bb{P}}\left(P\Phi_X^t(\theta_{p_0}),
      P\Phi_X^t(\theta^\prime)\right) > C,\quad\forall t\in\bb{R}, 
  \end{equation}
  where $\mathrm{dist}_{\bb{P}}$ denotes the distance function on
  $\bb{P}(NX)$ induced by the Riemannian structure
  $\langle\cdot,\cdot\rangle_{NX}$.

  Then, taking into account conservativeness proved in
  Lemma~\ref{lem:conserv-normal-flow},
  estimate~(\ref{eq:theta_p-theta-0-dist}) and
  equation~(\ref{eq:Y-hat-growth}), we have that any vector
  $\hat{v}\in NX_{p_0}$ whose $\mathrm{pr}_{X}$-projection is equal to
  $\theta^\prime\in K_0\cap\pi_{\bb{P}}^{-1}(p_0) $ will satisfies the
  following estimate:
  \begin{equation}
    \label{eq:hyperb-bundle-K-0}
    \left\|N\Phi_x^{t}(\hat{v})\right\|_{NX}\leq C^\prime 
    \exp\bigg(-t \int_M L_\Sigma\Omega\bigg)\|\hat{v}\|_{NX}
    ,\quad\forall t\in\bb{R},
  \end{equation}
  and for some real constant $C^\prime > 0$, that just depends on
  constant $C$ of estimate~(\ref{eq:theta_p-theta-0-dist}).

  From equation~(\ref{eq:Y-hat-growth}) and
  estimate~(\ref{eq:hyperb-bundle-K-0}) (and supposing that $\int_M
  L_\Sigma\Omega\neq 0$), we clearly conclude that Oseldets splitting
  (see \cite{oseledets}) of the linear cocycle $\{N\Phi_X^t\}$ is not
  just measurable, but continuous and uniformly hyperbolic. This
  implies that $\{N\Phi_X^t\}$ is an Anosov cocycle, and by
  Proposition~\ref{pro:Anosov-normal-impl-Anosov}, we know that $X$
  must be Anosov, which is clearly impossible, since $\{\Phi_X^t\}$
  does not have any periodic orbit.

  Therefore, the absurd comes from our supposition that $\int_M
  L_\Sigma\Omega$ could be non-null. Finally,
  equation~(\ref{eq:Y-hat-growth}) let us assure that $\hat{Y}$ is a
  $N\Phi_X$-invariant section, as desired.
\end{proof}

Now, we are ready for proving the theorem:

\begin{proof}[Proof of Theorem~\ref{thm:minimal-set-P(NX)}]
  Let $K_0\subset\bb{P}(NX)$, $p_0\in M$ and $\theta^\prime\in
  K_0\cap\pi_{\bb{P}}^{-1}(p_0)\in\bb{P}(NX)$ as above. Let
  $\hat{v}\in NX_{p_0}$ verifying
  $\mathrm{pr}_{\bb{P}}(\hat{v})=\theta^\prime$.

  We can rewrite estimate~(\ref{eq:theta_p-theta-0-dist}) as
  \begin{equation}
    \label{eq:Y-hat-v-hat-angle}
    \inf_{t\in\bb{R}} \sphericalangle\left(Y(\Phi_X^t(p)),
      N\Phi_X^t(\hat{v})\right) > 0.
  \end{equation}
  
  Putting together equation~(\ref{eq:invariance-hat-Y}),
  estimate~(\ref{eq:Y-hat-v-hat-angle}) and
  Lemma~\ref{lem:conserv-normal-flow}, we get that there exists a real
  constant $C^{\prime\prime}>1$ verifying
  \begin{equation}
    \label{eq:v-hat-bounded}
    \frac{1}{C^{\prime\prime}} < \left\|N\Phi_X^t(\hat{v})\right\|_{NX}
    < C^{\prime\prime},\quad\forall t\in\bb{R}.
  \end{equation}

  Now, consider another vector $\hat{w}\in NX_{p_0}\setminus\{0\}$
  such that $\mathrm{pr}_{\bb{P}}(\hat{w})\not\in K_\Sigma\cup K_0$.
  Since $K_\Sigma$ and $K_0$ are the only minimal sets for
  $\{P\Phi_X^t\}$, we know that the $\omega$-limit of
  $\mathrm{pr}_{\bb{P}}(\hat{w})$ must be either $K_0$ or $K_\Sigma$.
  Let us suppose that the positive semi-orbit of
  $\mathrm{pr}_{\bb{P}}(\hat{w})$ accumulates on $K_\Sigma$. This
  implies that
  \begin{equation}
    \label{eq:Y-hat-w-hat-angle}
    \lim_{t\to+\infty}\sphericalangle\left(
      \hat{Y}\big(\Phi_X^t(p_0)\big), N\Phi_X^t(\hat{w})\right) = 0.
  \end{equation}

  Once again, taking into account that $\{N\Phi_X^t\}$ preserves the
  symplectic form $\kappa$ and section $\hat{Y}\in\Gamma(NX)$, we see
  that equation~(\ref{eq:Y-hat-w-hat-angle}) implies that
  \begin{equation}
    \label{eq:w-hat-unbounded}
    \left\|N\Phi_x^t(\hat{w})\right\|\longrightarrow \infty,\quad
    \text{as } t\to+\infty.
  \end{equation}

  Finally, we clearly see that estimates~(\ref{eq:Y-hat-v-hat-angle}),
  (\ref{eq:v-hat-bounded}) and (\ref{eq:w-hat-unbounded}) violate
  conservativeness.

  We can analogously get a contradiction supposing that the
  $\omega$-limit of $\mathrm{pr}_{\bb{P}}(\hat{w})$ is $K_0$, and the
  we conclude that $K_\Sigma$ is the only minimal set for
  $\{P\Phi_X^t\}$.
\end{proof}

\subsection{Dynamics of the Normal Flow II}
\label{sec:normal-flow-II}

In paragraph~\ref{sec:normal-flow} we begun the analysis of the
dynamics of the normal flow $\{N\Phi_x^t\}$. After what we have just
done in paragraph~\ref{sec:proj-flow}, here we shall see that some of
the results previously gotten can be considerably improved. In fact,
we will completely characterize the dynamics of $\{N\Phi_X^t\}$,
showing that it exhibits a parabolic behavior.

In Lemma~\ref{lem:quasi-anosov} we showed that there was some non-null
vector in $NX$ such that its whole $N\Phi_X$-orbit was bounded. On the
other hand, in Lemma~\ref{lem:dynamics-on-E-Sigma}, under the
assumption that there were two different minimal sets for
$\{P\Phi_X^t\}$, we proved that there existed
$\hat{Y}\in\Gamma(E_\Sigma)$ which was invariant under the action of
$\{N\Phi_X^t\}$. Our first result of this paragraph consists in
proving that we can get the same invariant section assuming in this
case that $K_\Sigma$ is the only minimal set:

\begin{lemma}
  \label{lem:Y-hat-invariant}
  There exists a non-vanishing section $\hat{Y}\in\Gamma(E_\Sigma)$
  verifying
  \begin{equation}
    \label{eq:Y-hat-invariant}
    N\Phi_X^t\big(\hat{Y}(p)\big)=\hat{Y}(\Phi_X^t(p)),
    \quad\forall p\in M,\ \forall t\in\bb{R}.
  \end{equation}
\end{lemma}

\begin{proof}
  Let $\hat{Y}_0,\hat{Y}\in\Gamma(E_\Sigma)$ and $L_\Sigma\in
  C^\infty(M,\bb{R})$ be as in Lemma~\ref{lem:dynamics-on-E-Sigma}.

  Recalling equation~(\ref{eq:Y-hat-growth}), we have
  \begin{displaymath}
    N\Phi_X^t(\hat{Y}(p))=\exp\bigg(t\int_M L_\Sigma\Omega\bigg)
    \hat{Y}(\Phi_X^t(p)),\quad\forall t\in\bb{R}.
  \end{displaymath}

  On the other hand, by Lemma~\ref{lem:quasi-anosov}, we know that
  there exists $\hat{v}_0\in NX\setminus\{0\}$ which $N\Phi_X$-orbit
  is bounded, and applying Theorem~\ref{thm:minimal-set-P(NX)} we get
  \begin{equation}
    \label{eq:Y-hat-v-hat-0-angle}
    \lim_{t\to\pm\infty} \mathrm{dist}_{\bb{P}}\left(
      \mathrm{pr}_{\bb{P}}
      \bigg(\hat{Y}\Big(\Phi_X^t\big(\pi_{N}(\hat{v}_0) 
      \big)\Big)\bigg),
      \mathrm{pr}_{\bb{P}}\left(N\Phi_x^t(\hat{v}_0)\right)\right) =
    0. 
  \end{equation}

  This clearly implies that $\|N\Phi_X^t(\hat{Y})\|_{NX}$ cannot
  exhibit exponential growth, and therefore, $\int_M
  L_\Sigma\Omega=0$, getting the desired invariance of $\hat{Y}$.
\end{proof}

Next, notice that $\pi_N\colon NX\to M$ is an orientable vector bundle
with $2$-dimensional fibers and $\hat{Y}$ is non-singular section of
this bundle. This clearly implies that $\pi_N\colon NX\to M$ is
globally trivial, in particular, we can find a smooth section
$\hat{Z}_0\in\Gamma(NX)$ verifying
\begin{equation}
  \label{eq:hat-Z-0-definition}
  \kappa\left(\hat{Y}(p),\hat{Z}_0(p)\right)=1,\quad\forall p\in M,
\end{equation}
and in particular, it holds $\mathrm{span}\{\hat{Y},\hat{Z}_0\}=NX$.

Theorem~\ref{thm:minimal-set-P(NX)} let us affirm that there exists
$\sigma\in\{-1,1\}$ satisfying
\begin{equation}
  \label{eq:hat-Z-0-assymp-angle}
  \begin{split}
    \lim_{t\to+\infty}
    \sphericalangle\left(N\Phi_X^t\big(\hat{Z}_0(p)\big),
      \sigma\hat{Y}\big(\Phi_X^t(p)\big)\right) &=0, \\
    \lim_{t\to-\infty}
    \sphericalangle\left(N\Phi_X^t\big(\hat{Z}_0(p)\big),
      -\sigma\hat{Y}\big(\Phi_X^t(p)\big)\right) &=0.
  \end{split}
\end{equation}

There is no lost of generality if we suppose that $\sigma=1$ in
(\ref{eq:hat-Z-0-assymp-angle}).

Using $\{\hat{Y},\hat{Z}_0\}$ as an ordered basis for $NX$,
$N\Phi_X^t\colon NX_p\to NX_{\Phi_X^t(p)}$ can be represented as an
element of $\SL(2,\bb{R})$, and indeed, it will have the following
form:
\begin{equation}
  \label{eq:NPhi-in-hat-Y-hat-Z-0-basis}
  N\Phi_X^t(p)=
  \begin{pmatrix}
    1 & \hat{a}(p,t) \\
    0 & 1
  \end{pmatrix},
\end{equation}
where $\hat{a}\colon M\times\bb{R}\to\bb{R}$ is a smooth function
satisfying $\hat{a}(\cdot,0)=0$.

Then, if we define $\hat{A}\in C^\infty(M,\bb{R})$ by
$\hat{A}(p)\doteq\partial_t\hat{a}(p,t)\big|_{t=0}$, we can find a
smooth real function $\hat{B}$ verifying
\begin{equation}
  \label{eq:hat-B-def}
  \Lie_X\hat{B}= -\hat{A} + \int_M \hat{A}\Omega. 
\end{equation}

Function $\hat{B}$ can be used for defining a new section
\begin{displaymath}
  \hat{Z}\doteq \hat{Z}_0+\hat{B}\hat{Y}\in\Gamma(NX),
\end{displaymath}
and in this way we clearly have
\begin{equation}
  \label{eq:hat-Z-evolution}
  N\Phi_X^t(\hat{Z}(p))=\hat{Z}\big(\Phi_X^t(p)\big)+
  t\left(\int_M\hat{A}\Omega\right)
  \hat{Y}\big(\Phi_X^t(p)\big).
\end{equation}
for any $t\in\bb{R}$ and $p\in M$.

From (\ref{eq:hat-Z-0-assymp-angle}) and (\ref{eq:hat-Z-evolution}) we
easily see that $\int_M\hat{A}\Omega>0$, proving that in fact,
$\{N\Phi_X^t\}$ exhibits a parabolic behavior as desired.

\subsection{Dynamics on $\Sigma$}
\label{sec:dyn-on-Sigma}

In this short paragraph we shall analyze the dynamics of the flow
$D\Phi_X\colon TM\times\bb{R}\to TM$ restricted to the invariant
sub-bundle $\Sigma\to M$.

Our main result consists in proving that $\{D\Phi_X^t\}$ on
$\Sigma\subset TM$, as $\{N\Phi_X^t\}$ on $NX$, has a parabolic
behavior. In fact, the techniques used in here are very similar to
those used in paragraph~\ref{sec:normal-flow}. The only novelty is
that \emph{a priori} we do not have any information about the
projective flow induced by $D\Phi_X\colon\Sigma\times\bb{R}\to\Sigma$.

In this case we know that, for each $p$ and $t$,
$D\Phi_X^t(X(p))=X(\Phi_X^t(p))$ and therefore, we should prove that
all the vectors non-collinear with $X$ have polynomial growth and
their directions converge to the direction of $X$.

Let us start considering any smooth vector field
$Y_0\in\Gamma(\Sigma)\subset\X(M)$ verifying
\begin{equation}
  \label{eq:Y-0-pr-NX}
  \mathrm{pr}_X(Y_0(p))=\hat{Y}(p),\quad\forall p\in M.
\end{equation}

Then, notice that putting together
equations~(\ref{eq:invariance-hat-Y}) and (\ref{eq:Y-0-pr-NX}) we can
affirm that there exists a smooth function $A\in C^\infty(M,\bb{R})$
verifying
\begin{equation}
  \label{eq:Y-0-Lie-derivative}
  \Lie_X Y_0 = AX.
\end{equation}

Once again, since $X$ is cohomology-free, there exists $B\in
C^\infty(M,\bb{R})$ satisfying
\begin{equation}
  \label{eq:B-definition}
  \Lie_X B = - A + \int_M A\Omega.
\end{equation}

We use this function $B$ for defining a new vector field
\begin{equation}
  \label{eq:Y-definition}
  Y\doteq Y_0 + BX\in\Gamma(\Sigma)\subset\X(M).
\end{equation}

Notice that it continues to hold $\mathrm{span}\{X,Y\}=\Sigma\subset
TM$ and, additionally, we get
\begin{equation}
  \label{eq:Y-Lie-derivative}
  \Lie_X Y \equiv \left(\int_M A\Omega\right) X.
\end{equation}

Thus, we have the following
\begin{lemma}
  \label{eq:unboundedness-of-Y}
  Function $A\in C^\infty(M,\bb{R})$ given by
  equation~(\ref{eq:Y-0-Lie-derivative}) satisfies
  \begin{displaymath}
    \int_M A\Omega\neq 0.
  \end{displaymath}
\end{lemma}

\begin{proof}
  Contrarily, let us suppose that $\int_M A\Omega=0$.

  Then, equation~(\ref{eq:Y-Lie-derivative}) is equivalent to say that
  $[X,Y]\equiv 0$, \ie $X$ and $Y$ commute. Since $X$ and $Y$ generate
  $\Sigma$, in particular we have that they are everywhere linearly
  independent, and so, these vector fields induce a locally free
  $\bb{R}^2$-action on $M$.
  
  Finally, a classical result due to Harold Rosenberg, Robert
  Roussarie and David Weil~\cite{ros} affirms that the only orientable
  closed $3$-manifolds admitting locally free $\bb{R}^2$-actions are
  $2$-torus bundles over a circle, and our manifold $M$ clearly does
  not satisfies this property since we are assuming that
  $H_1(M,\bb{Q})=0$.
\end{proof}

As a corollary of this lemma we easily see that, given any $p\in M$,
it holds $\|D\Phi_X^t(Y(p))\|\to\infty$, uniformly as $t\to\pm\infty$,
and
\begin{equation}
  \begin{split}
    \label{eq:Y-convergence-angle}
    \lim_{t\to+\infty} \sphericalangle\left(D\Phi_X^{t}(Y(p)),\sigma_0
      X(\Phi_X^{t}(p))\right) &=0, \\
    \lim_{t\to-\infty}
    \sphericalangle\left(D\Phi_X^{t}(Y(p)),-\sigma_0
      X(\Phi_X^{t}(p))\right) &=0,
  \end{split}
\end{equation}
where $\sigma_0\doteq \sign\left(\int_M A\Omega\right)\in\{1,-1\}$ and
$\sphericalangle(\cdot,\cdot)$ stands for the angle (measured with
respect to the Riemannian structure $\langle\cdot,\cdot\rangle$)
between two non-null tangent vectors.

For the sake of simplicity, and since we do not loose any generality,
we shall assume that $\int A\Omega >0$, and thus, $\sigma_0=1$.

Summarizing what we have just proved, $D\Phi_X^t\colon
\Sigma_p\to\Sigma_{\Phi_X^t(p)}$ is a parabolic linear map, and taking
the ordered set $\{X,Y\}$ as basis of $\Sigma\subset TM$, we can
represent it by
\begin{equation}
  \label{eq:DPhi-matrix-on-Sigma}
  D\Phi_X^t\Big|_{\Sigma} = 
  \begin{pmatrix}
    1 & t\left(\int_M A\Omega\right) \\
    0 & 1
  \end{pmatrix}.
\end{equation}

\subsection{Expansiveness}
\label{sec:expansivness}

Let us start this paragraph recalling the definition of
\emph{expansive flow} due to Rufus Bowen and Peter
Walters~\cite{bowen-walters}:
\begin{definition}
  \label{def:expansive-flow}
  Given a compact metric space $(K,d)$, a continuous flow $\Psi\colon
  K\times\bb{R}\to K$ is called \emph{expansive} if it satisfies the
  following property:

  For every $\epsilon>0$, there is a $\delta>0$ such that if there
  exists a pair of points $x,y\in K$ and a homeomorphism
  $h\colon\bb{R}\to\bb{R}$ with $h(0)=0$ verifying
  \begin{equation}
    \label{eq:expansivness-definition}
    d(\Psi^t(x),\Psi^{h(t)}(y))< \delta,\quad\forall t\in\bb{R},
  \end{equation}
  then $y=\Psi^\tau(x)$, for some $\tau\in(-\epsilon,\epsilon)$.

  Moreover, we shall say that $\Psi$ is \emph{positively expansive}
  (respec.  \emph{negatively expansive}) if above condition is
  satisfied replacing $\bb{R}$ by $(0,+\infty)$ (respec.
  $(-\infty,0)$) in equation~(\ref{eq:expansivness-definition}). More
  precisely, if it holds $y=\Psi^\tau(x)$, for some
  $\tau\in(-\epsilon,\epsilon)$, whenever
  \begin{displaymath}
    d(\Psi^t(x),\Psi^{h(t)}(y))<\delta,\quad\forall t\in(0,+\infty)\
    (\forall t\in(-\infty,0)).
  \end{displaymath}
\end{definition}

Our main goal now consists in proving that our flow $\{\Phi_X^t\}$ is
positively (and in fact also negatively) expansive.

For this, let us start observing that in
paragraph~\ref{sec:normal-flow-II} we have constructed a smooth
section $\hat{Z}\in\Gamma(NX)$ that verifies
equation~(\ref{eq:hat-Z-evolution}), where $\int_M\hat{A}\Omega\neq 0$
(in fact, we have supposed that this constant is positive). Then, if
$Z\in\X(M)$ is any smooth vector field verifying
$\mathrm{pr}_X(Z)=\hat{Z}$, we will clearly have that for every $p\in
M$,
\begin{equation}
  \label{eq:Z-norm-evolution}
  \left\|D\Phi_X^t(Z(p))\right\|\to\infty, \quad\text{when }
  t\to\pm\infty,
\end{equation}
being the convergence uniform.

On the other hand, equations (\ref{eq:hat-Z-0-assymp-angle}) and
(\ref{eq:Y-convergence-angle}) let us affirm that (modulo our sign
assumptions made there) for every $p$ it holds
\begin{equation}
  \label{eq:Z-ang-evolution}
  \sphericalangle\left(D\Phi_X^t\big(Z(p)\big),
    X\big(\Phi_X^t(p)\big)\right)\to 0,\quad \text{when } t\to
  +\infty, 
\end{equation}
being this convergence uniform, too.

Then, taking into account that $\{X,Y,Z\}$ is a global basis for $TM$,
jointly with equations (\ref{eq:unboundedness-of-Y}),
(\ref{eq:Y-convergence-angle}), (\ref{eq:Z-norm-evolution}) and
(\ref{eq:Z-ang-evolution}), we easily get

\begin{proposition}
  \label{pro:Phi-expansive}
  The flow $\{\Phi_X^t\}$ is positively expansive.
\end{proposition}

And then we are very close to the end of our proof. In fact, as we
will shortly see, there is no closed $3$-manifold supporting
positively expansive flows. The essential tool for getting this is the
work due to Miguel Paternain~\cite{patern} about the existence of
stable and unstable foliations for expansive flows on $3$-manifolds.

Let us briefly recall Paternain's results. For this we need to
introduce some additional notation. Let $K$ be any closed manifold,
$\mathrm{dist}\colon K\times K\to\bb{R}$ be any distance compatible
with the topology of $K$ and $\Psi\colon K\times\bb{R}\to\bb{K}$ be a
continuous expansive flow.

As usual, given any $x\in K$, we can define its \emph{stable} and
\emph{unstable sets} writing
\begin{displaymath}
  \begin{split}
    W^s(x,\Psi) & \doteq\left\{y\in K :
      d\big(\Psi^t(x),\Psi^t(y)\big)\to 0,\
      \text{as } t\to+\infty\right\}, \\
    W^u(x,\Psi) & \doteq\left\{y\in K :
      d\big(\Psi^{-t}(x),\Psi^{-t}(y)\big)\to 0,\ \text{as }
      t\to+\infty\right\},
  \end{split}
\end{displaymath}
respectively.

Thus, we can precisely state

\begin{theorem}[Paternain~\cite{patern}]
  \label{thm:paternain-inv-foliations}
  If $K$ is a closed $3$-manifold and $\Psi$ is an expansive flow on
  $K$, then there exists a finite set (maybe empty) of periodic orbits
  $\gamma_1,\gamma_2,\ldots,\gamma_n$ of $\Psi$ such that the
  partitions
  \begin{displaymath}
    \scr{F}^\sigma=\left\{W^\sigma(x,\Psi) : x\in
      M\setminus\bigcup_{i=1}^{n} \gamma_i\right\},\quad\text{for }
    \sigma=s,u, 
  \end{displaymath}
  are $C^0$ codimension-two foliations on $M\setminus\bigcup\gamma_i$.
\end{theorem}

In our particular case the flow $\{\Phi_X^t\}$ has no periodic orbit,
and hence, since we have proved that it is positively expansive, in
particular, it is expansive and then, this theorem lets us affirm
that, given any point $p\in M$, the set $W^s(p,\Phi_X)$ does not just
reduce to $\{p\}$. This clearly contradicts the fact that
$\{\Phi_X^t\}$ is positively expansive, and we finish our proof.


\chapter{Final Remarks and Problems}
\label{chap:remarks-problems}

\section{On Manifolds with $\beta_1(M)=0$}
\label{sec:remarks-beta1-0}

\subsection{$3$-manifolds and Weinstein Conjecture}
\label{sec:remarks-weinstein-conj}

As it was already explained in the introduction of this work, the main
goal behind Katok Conjecture is to understand all possible
(topological and analytical) obstructions than can appear when we look
for smooth solutions of cohomological equations.

In Chapter~\ref{chap:b_1>0} we analyzed the existence of
cohomology-free vector fields on $3$-manifolds with non-zero first
Betti number. In all the stages of the proof of Theorem~A it was
rather clear how the topology of the manifold imposed different
obstructions for the existence of cohomology-free vector fields, and
all those obstructions let us completely characterize the supporting
manifold.

Unfortunately, the situation is not that clear when we have to prove
that there is no rational homological $3$-sphere supporting
cohomology-free vector fields. First, in
Section~\ref{sec:inv-one-form}, we proved that a hypothetical
cohomology-free vector field on such a manifold had to preserve a
non-singular $1$-form and the analysis of the existence of
obstructions was very satisfactory in the case that the kernel of the
invariant $1$-form was integrable
(Section~\ref{sec:comp-integrable-case}): solving some cohomological
equations we completely characterized the dynamics of the derivative
of the flow and we saw that there was no flow with that behavior on
the tangent bundle.

Nevertheless, when we had to analyze the case where the kernel of the
invariant $1$-form determined a contact structure, we just proved that
our hypothetical vector field was collinear with the Reeb vector field
induced by the invariant $1$-form, and then we finished our proof
invoking Taubes' work on Weinstein Conjecture. From a purely formal
point of view, this is a correct and complete proof, but if we take
into account the real goal behind Katok Conjecture, we cannot affirm
that this is a satisfactory one, because we are not understanding the
nature of the obstructions that appear in this case. This is mainly
due to the fact that Taubes' techniques used in \cite{taubes} are
extremely different to those used in the rest of this work.

Hence, it would be very desirable to complete the analysis that we
started in Section~\ref{sec:contact-struct-case} not invoking Taubes'
proof of Weinstein Conjecture, getting a more ``cohomological'' proof.

\subsection{Higher Dimensional Manifolds}
\label{sec:higher-dim}

As the reader could see in Chapter~\ref{chap:b_1>0},
Theorem~\ref{thm:fede-jana} due to Federico and Jana Rodr\'\i
guez-Hertz had a very important role in the proof of Theorem~A.

Nevertheless, if the first Betti number of our supporting manifold is
zero, then this result does not supply any non-trivial information.

Therefore, it seems reasonable to propose the following

\begin{problem}
  \label{prob:fiber-over-T1}
  Let $M$ be a closed $d$-manifold, with $d\geq 5$. Let us assume that
  there exists $X\in\X(M)$ cohomology-free.  Then, does there exist a
  good fibration for $X$ $p\colon M\to\bb{T}^1$? In particular, must
  it hold $\beta_1(M)\geq 1$?
\end{problem}

Another problem that seems to be very helpful (but difficult) for
understanding the dynamics of cohomology-free diffeomorphisms on
higher dimensional manifolds, is the following one proposed by Richard
Luz and Nathan dos Santos~\cite{luz-santos}:
\begin{problem}
  If $M$ is a closed manifold, $f\in\Diff(M)$ is cohomology-free and
  $n\in\bb{Z}\setminus\{0\}$, is it true that $f^n$ is
  cohomology-free?
\end{problem}

Motivated by this problem, we propose the following one for vector
fields:
\begin{problem}
  If $M$ is a closed manifold, $p\colon\tilde{M}\to M$ a $k$-fold
  covering with $k\geq 2$ and $X\in\X(M)$ cohomology-free, is the
  $p$-lift vector field $\tilde{X}\doteq p^*(X)\in\X(M)$
  cohomology-free?
\end{problem}

\section{Globally Hypoelliptic Vector Fields}
\label{sec:glob-hypo-vect}
 
In the theory of Partial Differential Equations there is a family of
smooth vector fields that has been extensively studied and that, a
priori, strictly contains the family of cohomology-free vector fields.
These are the \emph{globally hypoelliptic vector fields}:

\begin{definition}
  Let $M$ be a closed orientable manifold and $X\in\X(M)$. We say that
  $X$ is \emph{globally hypoelliptic} if given any $T\in\Dis(M)$, it
  holds
  \begin{displaymath}
    \Lie_X T \in C^\infty(M,\bb{R})\subset\Dis(M)\Rightarrow T\in
    C^\infty(M,\bb{R}). 
  \end{displaymath}
\end{definition}

It is very easy to see that every cohomology-free vector field is
indeed globally hypoelliptic, but a priori these two concepts are not
equivalent.

The first result concerning the classification of globally
hypoelliptic vector fields is due to Stephen Greenfield and Nolan
Wallach who proved in \cite{hypo-vect} that, modulo $C^\infty$
conjugacy, the constant vector fields on $\bb{T}^2$ verifying a
Diophantine condition like~(\ref{eq:diophantine-def}) are the only
examples on closed surfaces. This led them to propose the following

\begin{conjecture}[Greenfield-Wallach Conjecture~\cite{hypo-vect}]
  \label{conj:greenfield-wallach}
  Tori are the only closed manifolds that support globally
  hypoelliptic vector fields.
\end{conjecture}

It is interesting to remark that this implies Katok Conjecture. In
fact, Chen Wenyi and M. Y. Chi have proved in \cite{cc} that the only
globally hypoelliptic vector fields on tori are those smoothly
conjugated to constant vector fields satisfying a Diophantine
condition like~(\ref{eq:diophantine-def}).

However, one of the main results in \cite{cc} is Theorem 2.2 which
asserts that any globally hypoelliptic vector field on $\bb{T}^d$ is
cohomology-free. As Federico Rodr\'\i guez-Hertz has recently
observed, the proof of this result, presented by Chen and Chi in
\cite{cc}, continues to hold on any closed manifold, and consequently,
both families of vector fields coincide. Therefore, we have that
Greenfield-Wallach Conjecture and Katok Conjecture are indeed
equivalent.

\section{Positively Expansive Flows}
\label{sec:expansive-flows}

Given a compact metric space $(K,d)$ and a homeomorphism $h\colon K\to
K$, we say that $h$ is \emph{expansive} if there exists
$\varepsilon>0$ such that, for any pair of distinct points $x,y\in K$,
it holds
\begin{displaymath}
  \sup_{n\in\bb{Z}}d\left(f^n(x),f^n(y)\right)>\varepsilon,
\end{displaymath}
and we say that $h$ is \emph{positively expansive} if it holds
\begin{displaymath}
  \sup_{n\in\bb{N}_0}d\left(f^n(x),f^n(y)\right)>\varepsilon,
\end{displaymath}

It is a very well known fact that $h\colon K\to K$ is positively
expansive if and only if $K$ is a finite set.

On the other hand, in Section~\ref{sec:expansivness}, invoking a
result due to Miguel Paternain~\cite{patern}, we easily proved that
there does not exist any positively expansive flow on closed
$3$-manifolds.  However, we do not have any knowledge about the
existence of positively expansive flows on higher dimensional
manifolds. In fact, taking into account the simple classification of
positively expansive homeomorphisms, it seems natural to ask:
\begin{problem}
  If $\Psi\colon K\times\bb{R}\to K$ is a fixed-point free positively
  expansive flow, is it true that $K$ is homeomorphic to a finite
  disjoint union of copies of $\bb{T}^1$?
\end{problem}

\bibliographystyle{amsalpha}

\bibliography{biblio}
\end{document}